\documentclass[12pt, reqno]{amsart}
\usepackage{latexsym}
\usepackage{mathtools}
\usepackage{amsmath}
\usepackage{amsthm}
\usepackage{mathrsfs}
\usepackage{lmodern}
\usepackage[T1]{fontenc}
\usepackage{textcomp}
\usepackage[utf8]{inputenc}
\usepackage{color}
\usepackage{amssymb}
\usepackage{enumerate}
\usepackage[abbrev]{amsrefs}
\usepackage{indentfirst}
\usepackage{graphicx}
\usepackage{bmpsize}
\usepackage{hyperref}
\usepackage{enumitem}
\usepackage{enumerate}
\usepackage{url}
\usepackage{autobreak}
\usepackage[justification=justified]{caption}
\usepackage{bm}
\usepackage{here}
\usepackage[labelformat=empty,subrefformat=parens]{subcaption}
\renewcommand{\thefootnote}{} 

\theoremstyle{plain} 
\newtheorem{theorem}{\indent\sc Theorem}[section]
\newtheorem{lemma}[theorem]{\indent\sc Lemma}
\newtheorem{corollary}[theorem]{\indent\sc Corollary}
\newtheorem{proposition}[theorem]{\indent\sc Proposition}

\theoremstyle{definition} 
\newtheorem{definition}[theorem]{\indent\sc Definition}
\newtheorem{remark}[theorem]{\indent\sc Remark}
\newtheorem{example}[theorem]{\indent\sc Example}

\newenvironment{sproof}{%
  \proof}{\endproof}
\makeatletter
  
  \@addtoreset{equation}{section}
\makeatother
\newcommand{\cP}{\mathcal{P}}
\newcommand{\cD}{\mathcal{D}}
\newcommand{\cR}{\mathcal{R}}
\newcommand{\cS}{\mathcal{S}}
\newcommand{\bR}{\mathbb{R}}
\newcommand{\bS}{\mathbb{S}}

\newcommand\cinput[2]{\lower#1pt\hbox{\input{#2}}}

\makeatletter
\@namedef{subjclassname@2020}{\textup{2020} Mathematics Subject Classification}
\makeatother

\begin{document}

\title[Virtual Thompson's group]{VIRTUAL THOMPSON'S GROUP}
\author{Yuya Kodama and Akihiro Takano}
\date{}
\renewcommand{\thefootnote}{\arabic{footnote}}  
\setcounter{footnote}{0} 
\thanks{The first author was supported by JST, the establishment of university
fellowships towards the creation of science technology innovation, Grant Number JPMJFS2139.}
\begin{abstract}
For virtual knot theory, the virtual braid group was defined by generalizing the braid group. 
It was proved that any virtual link can be obtained by the closure of a virtual braid. 
On the other hand, due to work by Jones et al., it is known that any (oriented) link is constructed from an element of Thompson's group $F$. 
In this paper, we define the ``virtual version'' of Thompson's group $F$ and prove that any virtual link is constructed from an element of the group. 
\end{abstract}
\maketitle
\section{Introduction}
Virtual knot theory, introduced by Kauffman \cite{kauffman1999virtual}, is a generalization of classical knot theory.
There are some motivations in this theory.
One is knot theory in $\Sigma_g \times [0,1]$, where $\Sigma_g$ is a closed oriented surface of genus $g \geq 1$.
Classical knot theory can be regarded as the case of $g = 0$.
Another is complete correspondence with the Gauss diagrams, which are used to define a finite type invariant \cite{goussarov2000virtual}.
As in the braid group in classical knot theory, the virtual braid group is defined and studied.
Kamada \cite{kamada2007virtual}, and Kauffman and Lambropoulou \cite{kauffman2004virtual} introduced this notion and proved Alexander's theorem, that is, any virtual link can be obtained from the closure of a virtual braid.
Moreover, they showed Markov's theorem.
In other words, this theorem gives the necessary and sufficient condition that two braids have equivalent closures.

Recently, Jones \cite{jones2017thompson} introduced a method of constructing a link from an element of Thompson's group $F$, and proved Alexander's theorem.
It means that any link can be obtained from an element of $F$.
For the oriented case, Jones defined a subgroup $\overrightarrow{F}$ of $F$ whose element yields an oriented link and showed the theorem with the weaker version. 
After that Aiello \cite{aiello2020alexander} proved it completely.
Golan and Sapir \cite{golan2017jones} showed the subgroup $\overrightarrow{F}$ is isomorphic to $3$-adic Thompson's group $F(3)$.

Thompson's group $F$ is defined by Richard Thompson in $1965$. 
This group is known to be related to various areas and has been studied using various definitions such as piecewise linear maps on $[0, 1]$, pairs of binary trees, and so on. 
We consider $F$ as a diagram group by referring to the approach in \cite{golan2017jones}. 
The notion of diagram groups was suggested by Meakin and Sapir, and then Kilibarda \cite{kilibarda1994algebra} studied the groups for the first time. 
This class of groups has been well studied not only algebraically but also geometrically. 
For instance, these groups are finitely presented \cite{guba1997diagram}, torsion-free \cite{guba1997diagram}, totally orderable \cite{guba2006diagram}, and act freely and cellularly on a CAT(0) cubical complex \cite{farley2003finiteness}. 

In this paper, we generalize Thompson's group $F$ from the viewpoint of virtual knot theory. 
Namely, we define virtual Thompson's group $VF$ as a diagram group and show the following: 
\begin{theorem}\label{main_theorem}
Any virtual link can be obtained from an element in $VF$.
\end{theorem}

This paper is organized as follows: 
In Section \ref{section_preliminaries}, we first summarize definitions of virtual links, diagram groups, and Thompson's group $F$. 
Then we define virtual Thompson's group $VF$ as a diagram group. 
At the end of this section, we discuss some properties of diagram groups, and hence of $VF$. 
In Section \ref{section_construction}, we introduce a method of constructing a virtual link from an element in $VF$. 
This method is a generalization of the one for $F$. 
Then we discuss the relationship between elements of $VF$ and labeled binary trees. 
Some elements of $VF$ are represented by labeled binary trees. 
In this sense, we can regard $VF$ as a generalization of $F$. 
In Section \ref{section_proof_main_theorem}, we show that any virtual link is obtained from some element in $VF$. 
Similar to \cite{jones2017thompson, aiello2020alexander}, this is achieved by constructing the Tait graph from a virtual link and deforming it. 

Various other generalizations of Thompson's group $F$ are also known \cite{brin2007algebra,dehornoy2006group,brin2004higher,lodha2016nonamenable}. 
It is an interesting problem to study the relationship between them and $VF$. 
\section{Preliminaries} \label{section_preliminaries}
\subsection{Virtual knots and links} \label{subsection_virtualintro}
In this section, we give a short description about the virtual links.

\begin{figure}[tbp]
\begin{tabular}{cc}
\begin{minipage}{0.33\hsize}
\begin{center}
\includegraphics[height=50pt]{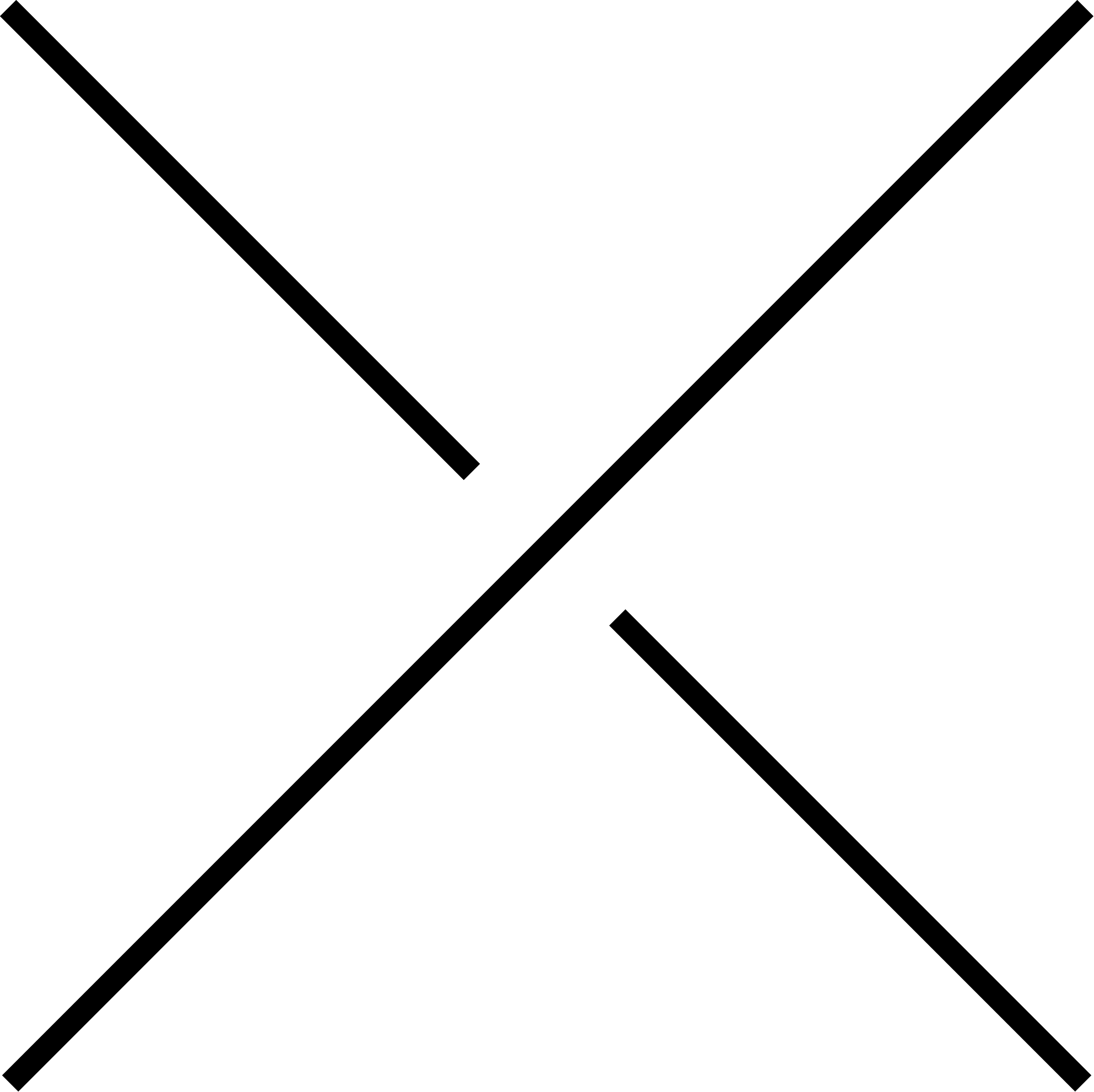}
\subcaption{Classical crossing}
\end{center}
\end{minipage}&
\hspace{-20pt}
\begin{minipage}{0.33\hsize}
\begin{center}
\includegraphics[height=50pt]{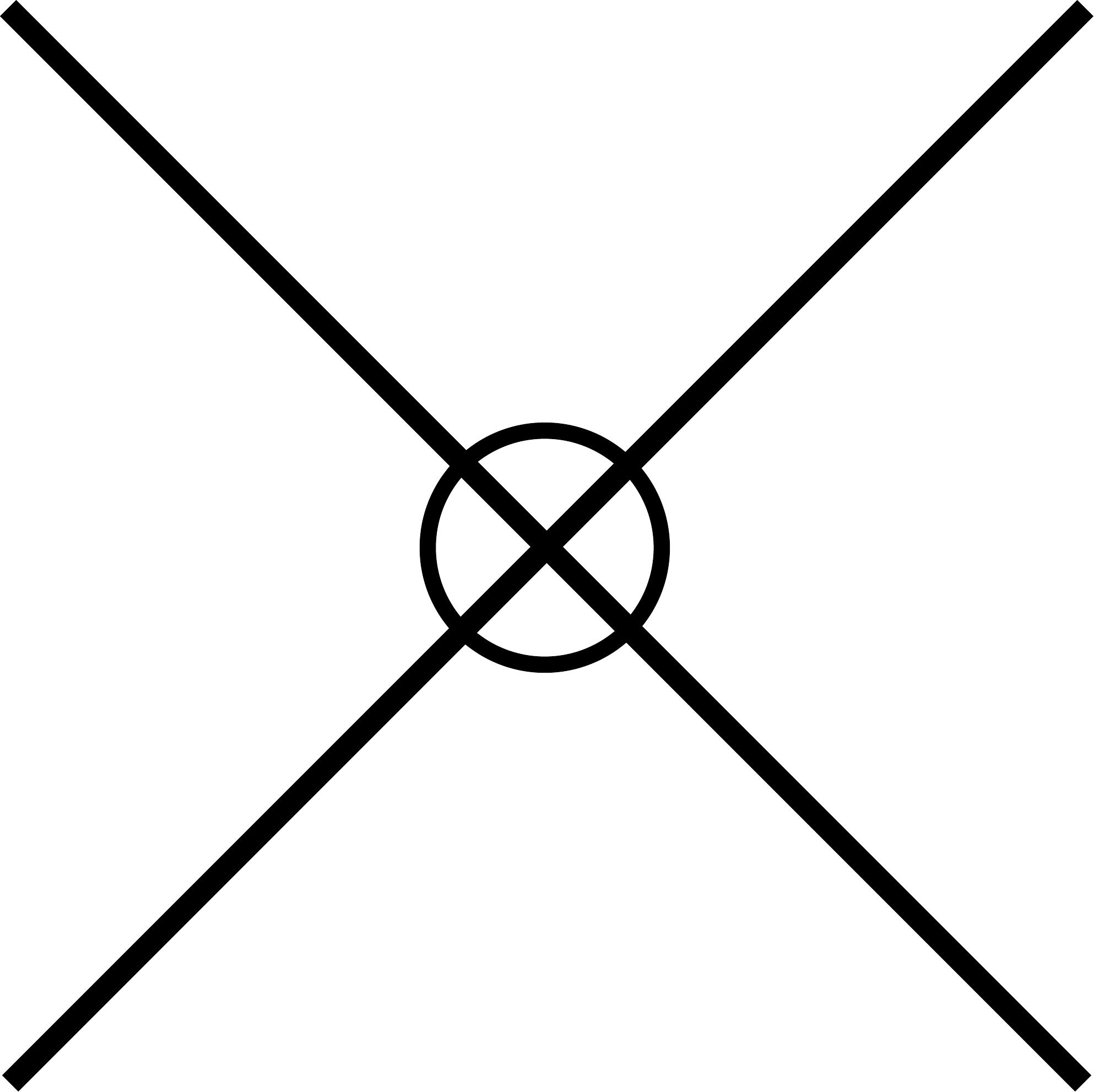}
\subcaption{Virtual crossing}
\end{center}
\end{minipage}
\end{tabular}
\vspace{-5pt}
\caption{Classical and virtual crossings}
\label{crossing}
\end{figure}

\begin{definition}
An $n$-component \textit{virtual link diagram} is an immersion of $n$ circles in $\bS^2 \left(= \bR^2 \cup \{ \infty \}\right)$ such that the multiple point set consists of finite number of transverse double points and each of them is labeled, either as a \textit{classical crossing} or as a \textit{virtual crossing} (see Figure \ref{crossing}).
In particular, if $n=1$, we also call it a \textit{virtual knot diagram}.
A virtual link diagram without virtual crossings is said to be \textit{classical}.
\end{definition}

\begin{figure}[tbp]
\begin{tabular}{ccc}
\hspace{-20pt}
\begin{minipage}{0.33\hsize}
\begin{center}
\includegraphics[height=70pt]{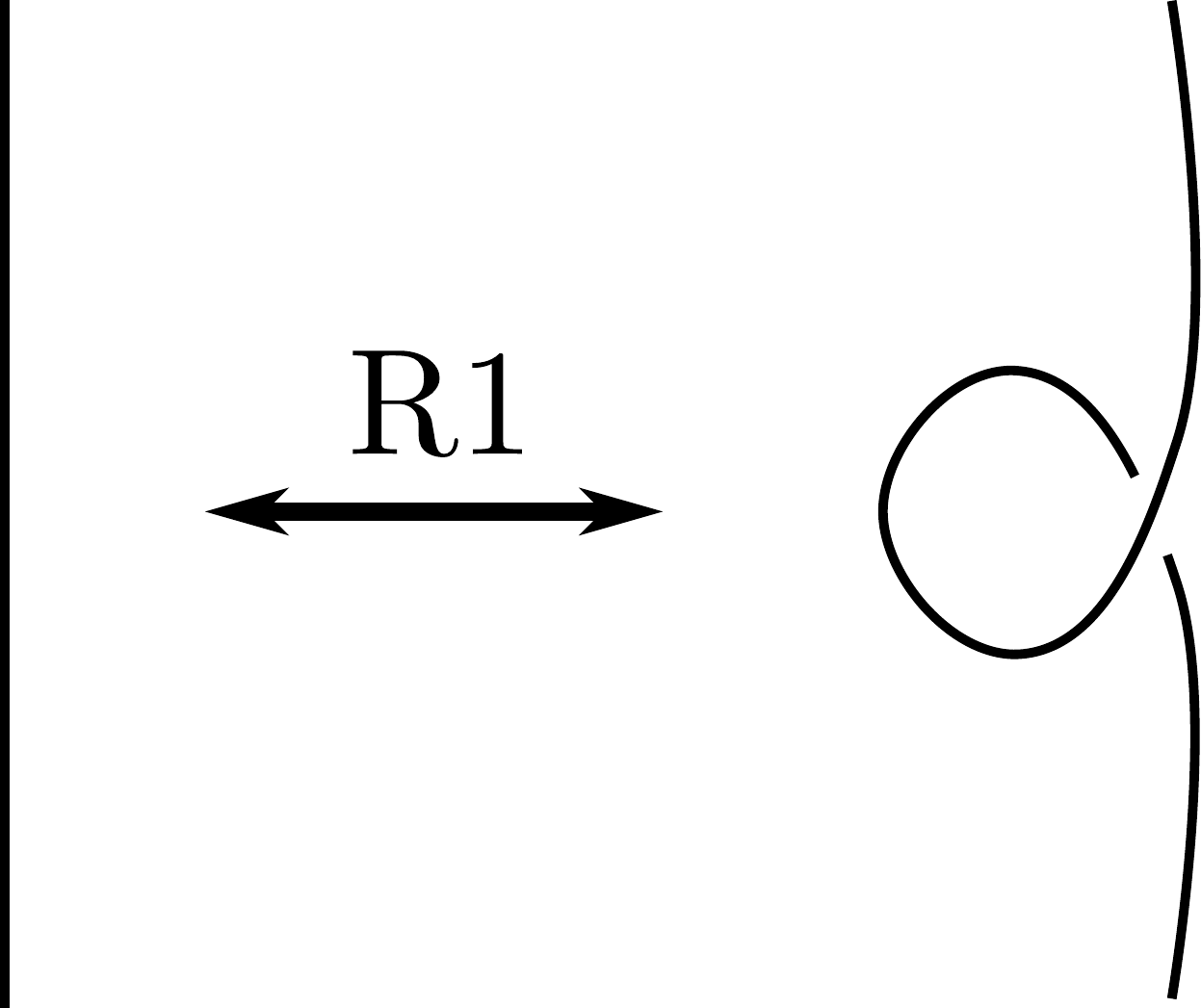}
\subcaption{}
\end{center}
\end{minipage}&
\hspace{-30pt}
\begin{minipage}{0.33\hsize}
\begin{center}
\includegraphics[height=70pt]{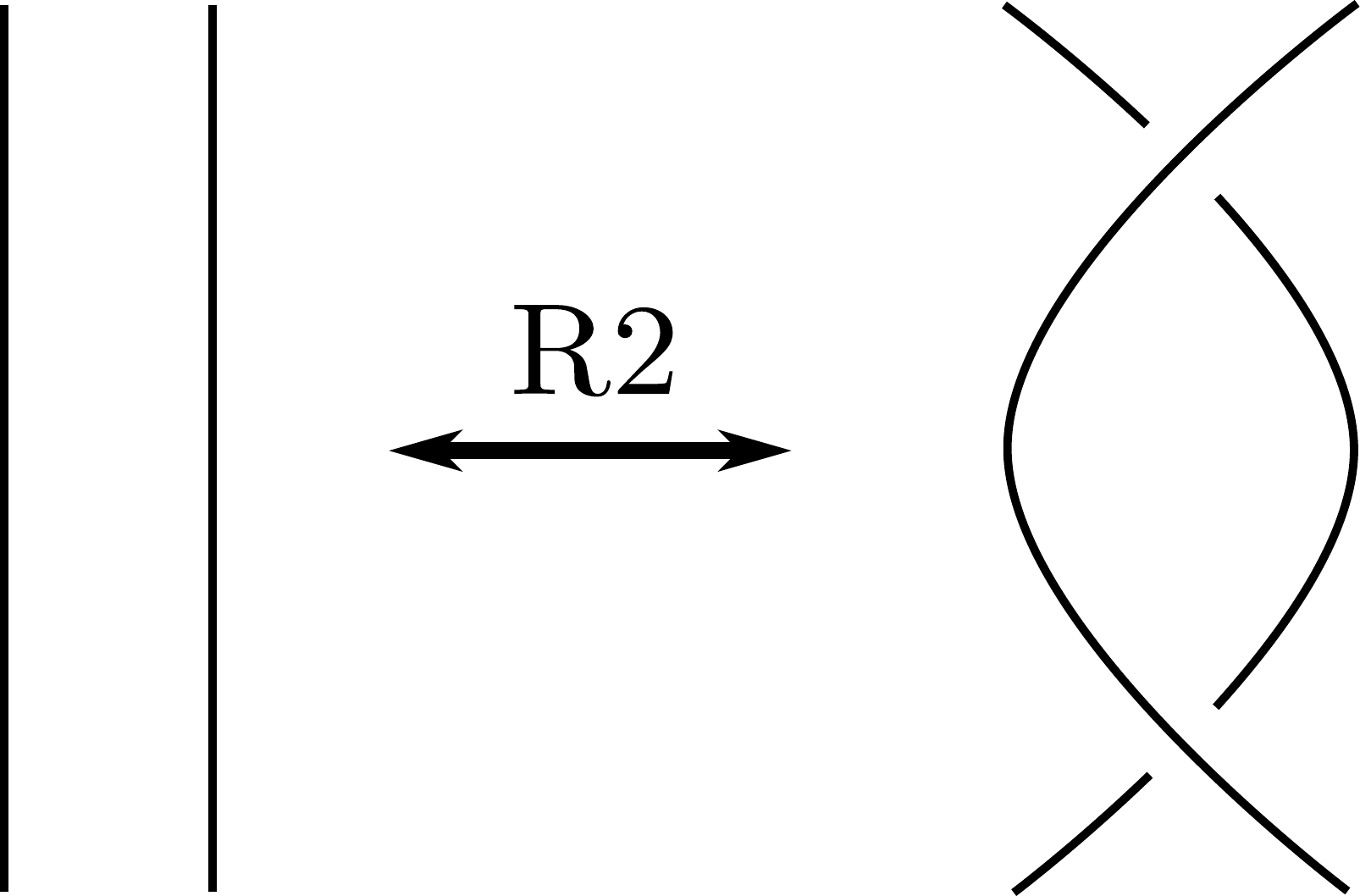}
\subcaption{}
\end{center}
\end{minipage}&
\begin{minipage}{0.33\hsize}
\begin{center}
\includegraphics[height=70pt]{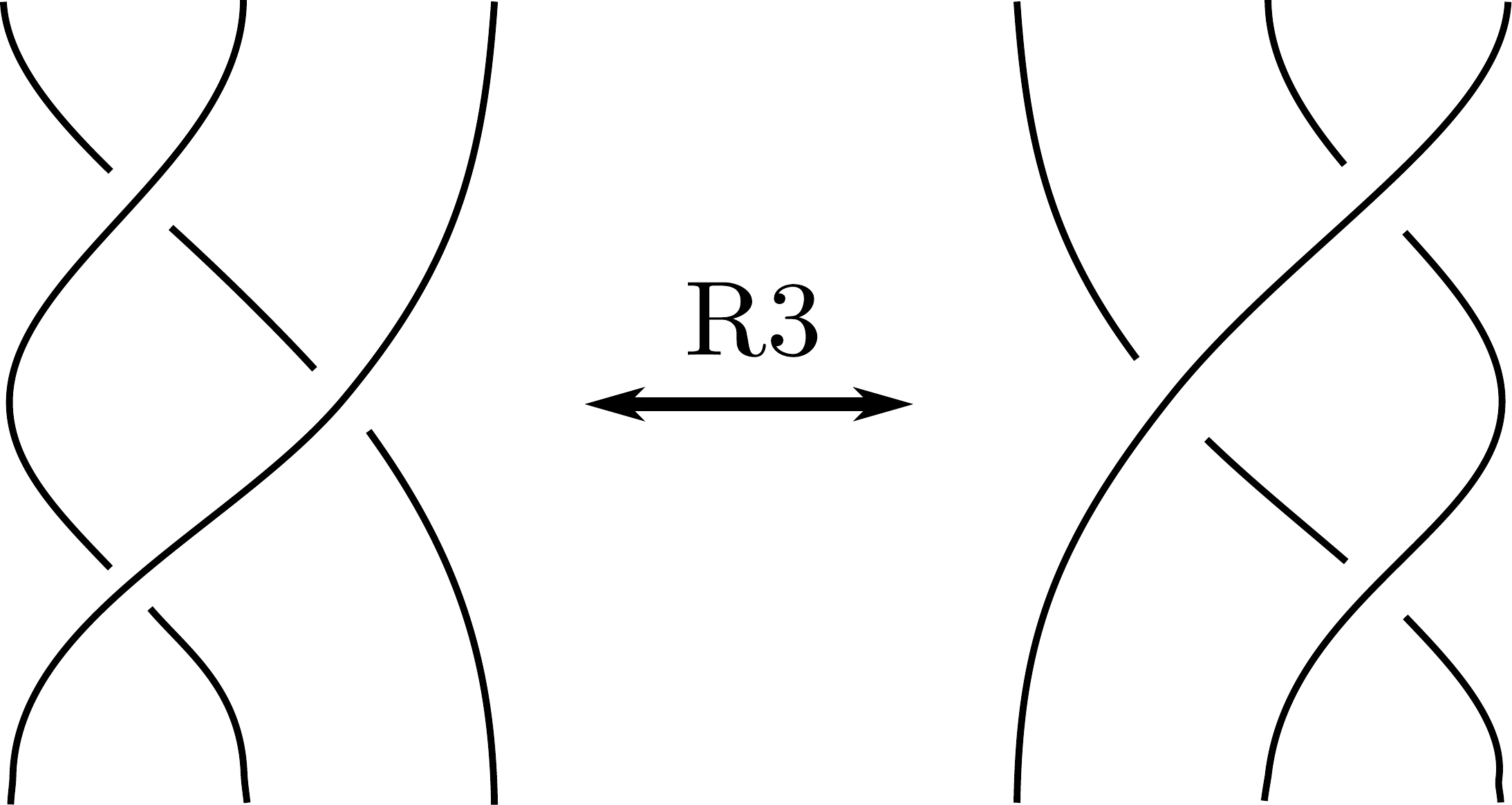}
\subcaption{}
\end{center}
\end{minipage}
\end{tabular}
\vspace{-25pt}
\caption{Classical Reidemeister moves}
\label{classical}
\end{figure}

\begin{figure}[tbp]
\begin{tabular}{ccc}
\hspace{-20pt}
\begin{minipage}{0.33\hsize}
\begin{center}
\includegraphics[height=70pt]{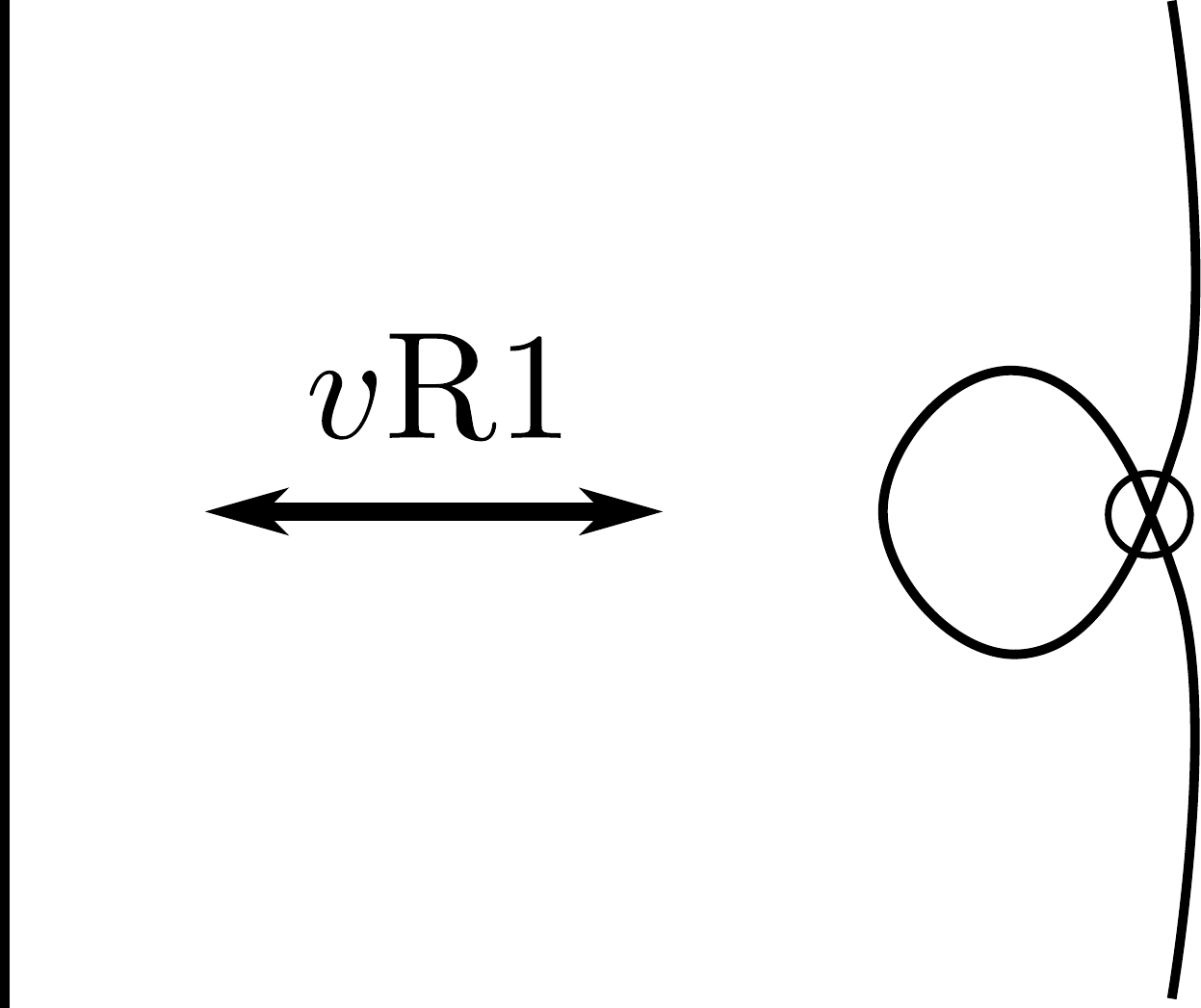}
\subcaption{}
\end{center}
\end{minipage}&
\hspace{-30pt}
\begin{minipage}{0.33\hsize}
\begin{center}
\includegraphics[height=70pt]{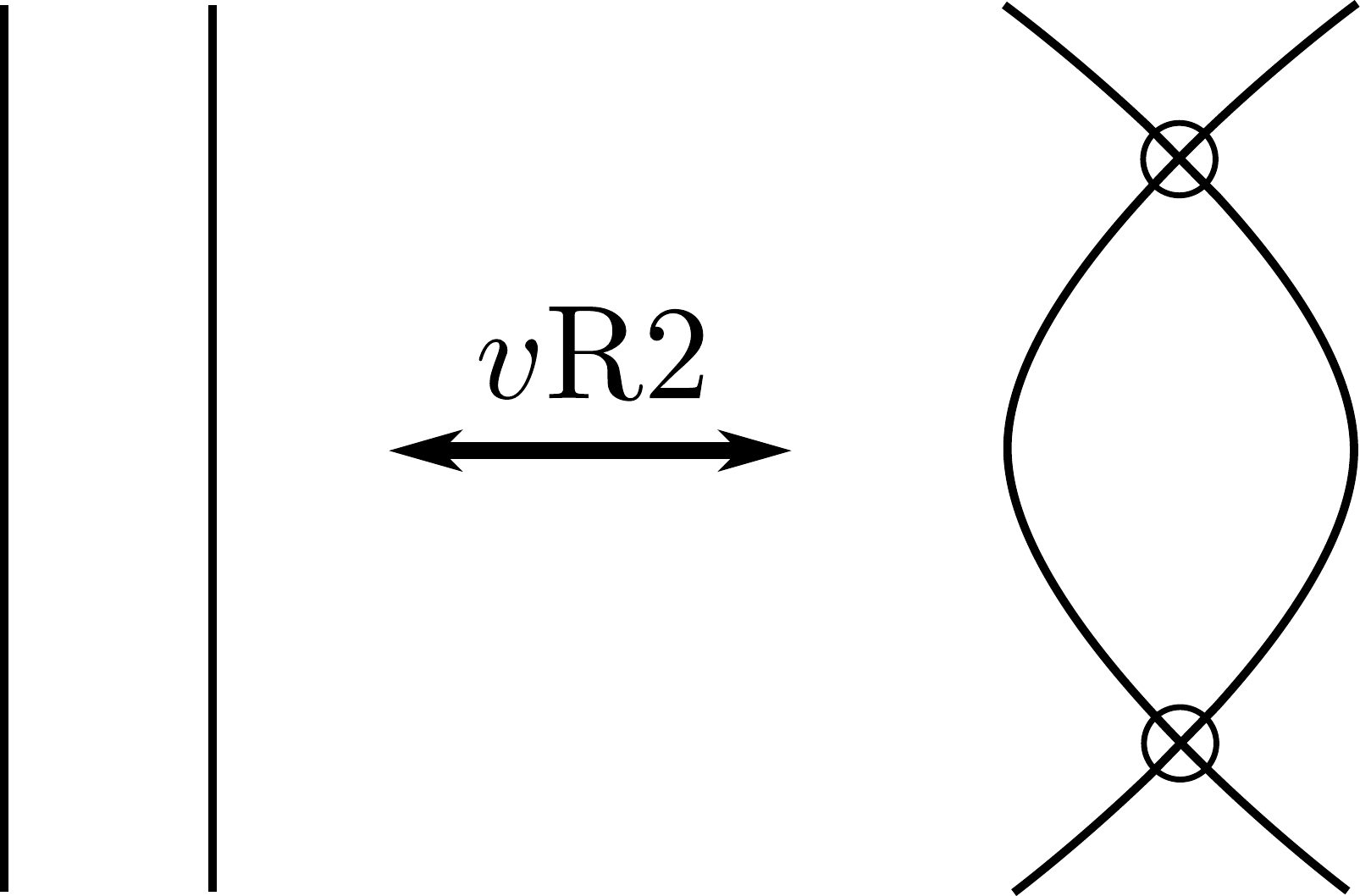}
\subcaption{}
\end{center}
\end{minipage}&
\begin{minipage}{0.33\hsize}
\begin{center}
\includegraphics[height=70pt]{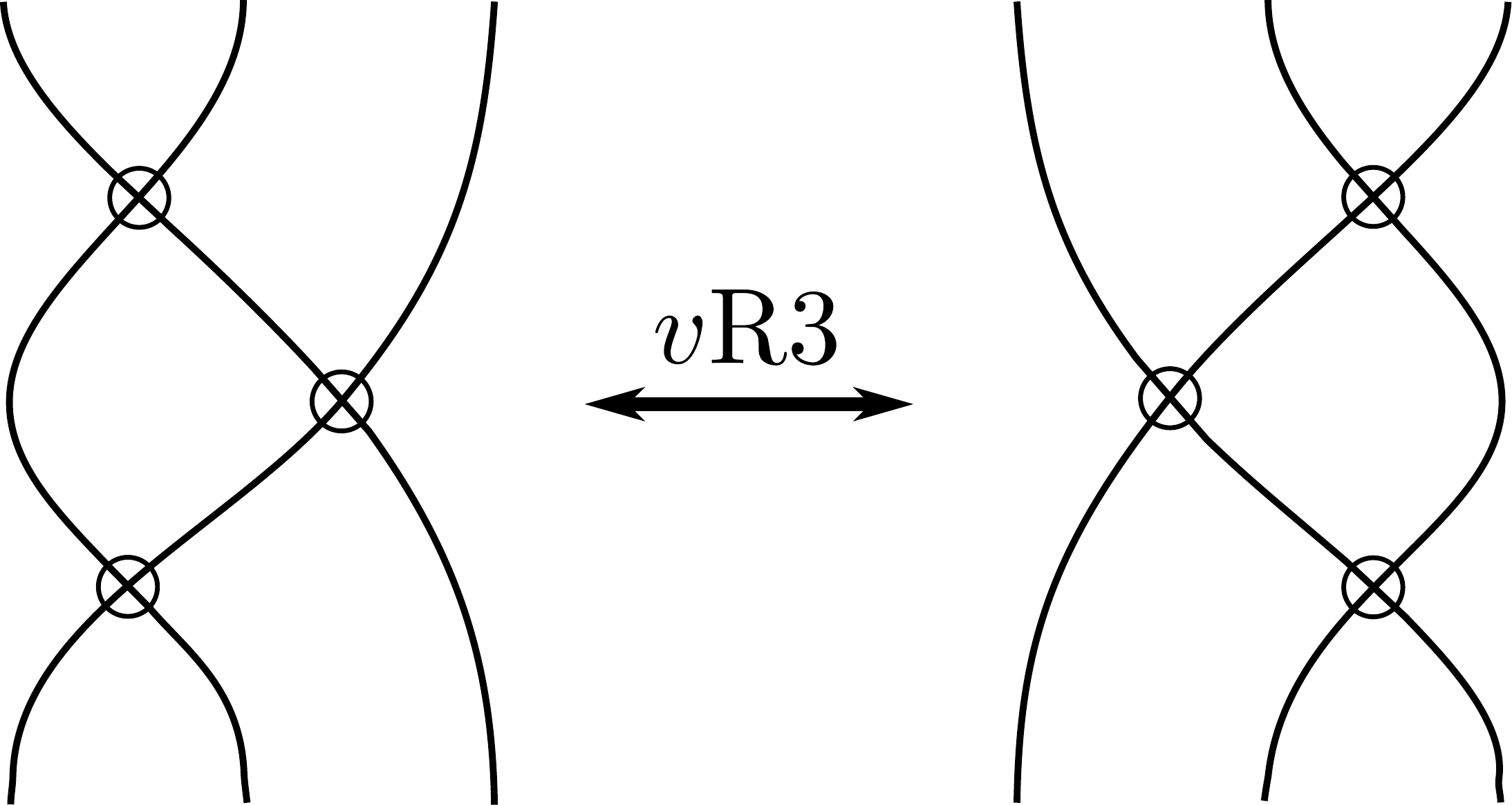}
\subcaption{}
\end{center}
\end{minipage}
\end{tabular}
\vspace{-25pt}
\caption{Virtual Reidemeister moves}
\label{virtual}
\end{figure}

\begin{figure}[tbp]
\begin{center}
\includegraphics[height=70pt]{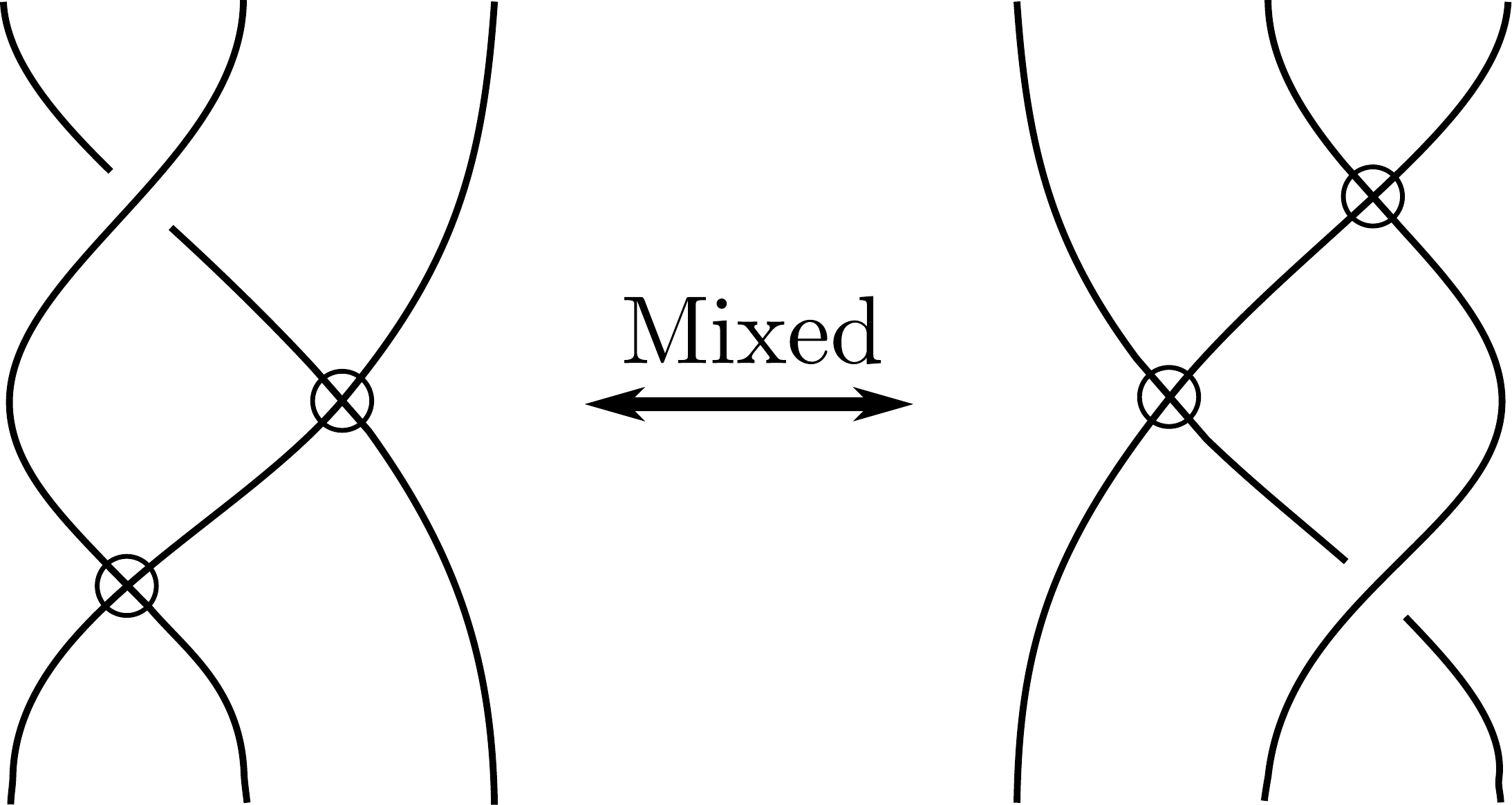}
\end{center}
\vspace{-5pt}
\caption{Mixed move}
\label{Mixed_OC}
\end{figure}

\begin{definition}
An $n$-component \textit{virtual link} is an equivalence class of the set of all $n$-component virtual link diagrams under the ambient isotopy on the plane and the generalized Reidemeister moves, that is, the (classical) Reidemeister moves (Figure \ref{classical}), the virtual Reidemeister moves (Figure \ref{virtual}), and the mixed move (Figure \ref{Mixed_OC}).
If $n=1$, we also call it a \textit{virtual knot}.
\end{definition}

Similarly to classical knot theory, it is natural to consider the notion of the virtual braid.
Kamada \cite{kamada2007virtual}, and Kauffman and Lambropoulou \cite{kauffman2004virtual} defined the virtual braid group and independently proved Alexander's theorem for virtual links:

\begin{theorem}[{\cite[Proposition 3]{kamada2007virtual}, \cite[Theorem 1]{kauffman2004virtual}}]
Any virtual link can be described as the closure of a virtual braid.
\end{theorem}
\subsection{Diagram groups over semigroups}
In this section, we briefly review the definition of diagram groups. 
Although our purpose is to define one diagram group, we explain the formal definition of diagram groups for the reader's convenience. 
See \cite{guba1997diagram} for details. 

Let $\cP= \langle \Sigma \mid \cR \rangle$ be a semigroup presentation. 
Here, $\Sigma$ is a finite set of generators, and $\cR$ is a finite set of relations of the form $u \to v$ where $u$ and $v$ are finite words on $\Sigma$. 
We always assume that there exists no relation of the form $u=u$, where $u$ is a finite word on $\Sigma$. 
For simplicity, if $u \to v$ is in $\cR$, then we regard $v \to u$ as also being in $\cR$. 

We fix a finite word $w$ on $\Sigma$. 
Roughly speaking, for the given word $w$, each element (called \textit{diagram}) of the diagram group represents a way of rewriting by relations from $w$ to itself again. 
Formally, for $w$, we define a diagram as a finite sequence of words on $\Sigma$ with the following form
\begin{align*}
w=w_1 \to w_2 \to \cdots \to w_{n-1} \to w_n=w
\end{align*}
where each $w_{i-1} \to w_i$ is of the form $w^\prime (p_{i-1})w^{\prime \prime} \to w^\prime (p_{i})w^{\prime \prime}$ and $p_{i-1} \to p_i$ is in $\cR$. 
We call each replacement of the word in the sequence \textit{cell} of the diagram. 

We define a reduction of a dipole as follows: 
Let $w=w_1 \to \cdots \to w_n=w$ be a diagram and assume that there exists $i$ such that $w_{i-1} \to w_{i} \to w_{i+1}$ is of the form $w^\prime (p_{i-1})w^{\prime \prime} \to w^\prime (p_{i})w^{\prime \prime} \to w^\prime (p_{i+1})w^{\prime \prime}$ and $p_{i-1}=p_{i+1}$ holds. 
In this case, we obtain a new diagram by eliminating $w_{i-1}$ and $w_{i}$, that is, by setting 
\begin{align*}
w=w_1 \to w_2 \to \cdots \to w_{i-2} \to w_{i+1} \to \cdots \to w_n=w. 
\end{align*}
This operation and its inverse are called \textit{reduction} and \textit{insertion of dipoles}, respectively. 
Combining these operations, we define the equivalence relation on the set of diagrams generated by the operations and the requirement that separated cells are commute. 
We call the group consisting of all equivalence classes of diagrams the \textit{diagram group} $\cD(\cP, w)$. 
The product on $\cD(\cP, w)$ is defined as follows: 
For two diagrams $w=a_1 \to \cdots \to a_n=w$ and $w=b_1 \to \cdots \to b_m=w$, we define their product to be the equivalence class of the ``concatenation''
\begin{align*}
w=a_1 \to \cdots \to a_n=w=b_1 \to \cdots \to b_m=w. 
\end{align*}
This product is well-defined. 

If we can not apply the reduction of dipoles, we call the diagram \textit{reduced}. 
For each element of $\cD(\cP, w)$, there exists a unique representative with reduced \cite{kilibarda1994algebra}. 
\begin{example} \label{example_diagram_rewritesystem}
Let $\cP=\langle a, b \mid a \to ab, b \to aa, a \to aa \rangle$ and $w=a$. 
Then
\begin{align*}
(a) \to (aa)=(a)a \to (aa)a=a(aa) \to a(a)=(aa) \to (a) 
\end{align*}
and
\begin{align}
(a) \to (ab)=a(b) \to a(aa) \to a(a)=(aa) \to (a) \label{reduction_dipole_word}
\end{align}
are reduced diagrams. 
Their product is
\begin{align*}
\begin{autobreak}
(a) 
\to (aa)=(a)a \to (aa)a=a(aa) \to a(a)=(aa) \to (a) 
\to (ab)=a(b) \to a(aa) \to a(a)=(aa) \to (a), 
\end{autobreak}
\end{align*}
and this diagram is also reduced. 
The diagram
\begin{align*}
(a) \to (ab)=(a)b \to (ab)b \to (a)b=a(b) \to a(aa) \to a(a)=(aa) \to (a)
\end{align*}
is not reduced since there exists $(a)b \to (ab)b \to (a)b$. 
If we reduce a dipole of this diagram, then we get diagram \ref{reduction_dipole_word}. 
The diagrams
\begin{align*}
\begin{autobreak}
(a)
 \to (ab)=a(b) \to a(aa)=(a)aa \to (aa)aa=aa(aa) 
\to aa(a)=a(aa) \to a(a)=(aa) \to (a)
\end{autobreak}
\end{align*}
and
\begin{align*}
\begin{autobreak}
(a)
 \to (ab)=(a)b \to (aa)b=aa(b) \to aa(aa) 
\to aa(a)=a(aa) \to a(a)=(aa) \to (a)
\end{autobreak}
\end{align*}
are equivalent since the cells $a(b) \to a(aa)$ and $(a)b \to (aa)b$ are separated. 
\end{example}
The notions in this section can also be represented by oriented graphs. 
Let $w=w_1w_2 \cdots w_n$ be a word where each $w_i$ is in $\Sigma$. 
We first define a trivial geometric diagram as follows: 

Let $v_1, v_2, \dots, v_{n+1}$ be a vertices, and each $v_i$ is connected to $v_{i+1}$ in this orientation. 
Namely, this graph consists of $n$ edges $(v_1, v_2), (v_2, v_3), \dots, (v_n, v_{n+1})$. 
We label each $(v_i, v_{i+1})$ as $w_i$ and omit the labels of vertices. 
We define this graph as a \textit{trivial geometric diagram} of $w$. 
For an oriented graph, given a vertex $v$ and two edges $(v^\prime, v)$ and $(v, v^{\prime \prime})$, we say that $(v^\prime, v)$ and $(v, v^{\prime \prime})$ are \textit{incoming} and \textit{outgoing edges} of $v$, respectively (cf.\ Figure \ref{incoming_outcoming}). 

Next, we define a geometric cell. 
Let $p \to q$ be in $\cR$ where $p=p_1 \cdots p_n$ and $q=q_1 \cdots q_m$. 
Let $v_{p_1}$ and $v_{p_n}$ be the vertices of the trivial diagram of $p$ such that the edges labeled by $p_1$ and $p_n$ are outgoing and incoming edges of $v_{p_1}$ and $v_{p_n}$, respectively.
For $q$, we define the vertices $v_{q_1}$ and $v_{q_m}$ similarly. 
Then the graph obtained by gluing $v_{p_1}$ and $v_{p_n}$ to $v_{q_1}$ and $v_{q_m}$, respectively, is called a \textit{geometric $(p, q)$-cell}. 

Generally, a geometric diagram is represented by attaching geometric cells to a trivial geometric diagram along corresponding sub-words successively. 
Let $w$ be a given word on $\Sigma$ and $w=w_1 \to w_2 \to \cdots w_{n-1} \to w_n=w$ be a diagram. 
Note that each $w_{i-1} \to w_i$ is of the form $w^\prime (p_{i-1})w^{\prime \prime} \to w^\prime (p_{i})w^{\prime \prime}$. 
Therefore we can attach a $(p_1, p_2)$-cell along the subgraph of a trivial geometric diagram of $w$ corresponding to $p_1$. 
Then the obtained graph has two paths, $w^\prime (p_{1})w^{\prime \prime}$ and $w^\prime (p_{2})w^{\prime \prime}$. 
Regarding $w^\prime (p_{2})w^{\prime \prime}$ as a trivial geometric diagram, and proceeding similarly to the end, we obtain a graph. 
We call this graph \textit{geometric diagram}. 
The trivial geometric (sub)diagrams corresponding to the top and bottom $w$ of the geometric diagram are called the \textit{top} and \textit{bottom paths}, respectively. 
In the top path or bottom path, the vertex with only one outgoing or incoming edge is called the \textit{initial} or \textit{terminal vertex}, respectively.
Note that for each diagram, there exists exactly one initial vertex and one terminal vertex since the top path and bottom path share them. 

\begin{figure}[tbp]
	\begin{center}
		\includegraphics[height=70pt]{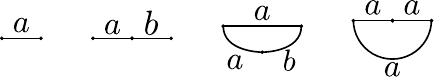}
	\end{center}
\caption{Trivial geometric diagrams of $a$ and $ab$, and geometric $(a, ab)$-cell and $(aa, a)$-cell}
\label{geo_triv_cell}
\end{figure}

\begin{figure}[tbp]
	\begin{center}
		\includegraphics[height=100pt]{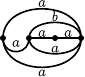}
	\end{center}
\caption{The geometric diagram corresponding to diagram \ref{reduction_dipole_word}}
\label{geo_diagram}
\end{figure}

Similar to the equivalence relation of diagrams, we define the equivalence relation on the set of geometric diagrams. 
In the rest of this paper, we do not distinguish between geometric diagrams and diagrams. 
See Figures \ref{geo_triv_cell}, \ref{geo_diagram} for examples of trivial geometric diagram, geometric cell, and geometric diagram corresponding to Example \ref{example_diagram_rewritesystem}. 
\subsection{Thompson's group $F$ as a diagram group and its generalization}
In this section, we recall Thompson's group $F$, then give the definition and some properties of virtual Thompson's group $VF$. 
This group is the most important one in this paper. 

We first outline the definition of Thompson's group $F$. 
It is known that there exist various (equivalent) definitions for this group. 
We define this group as pairs of binary trees and then see the correspondence of the other realizations. 

\begin{figure}[tbp]
	\begin{center}
		\includegraphics[height=80pt]{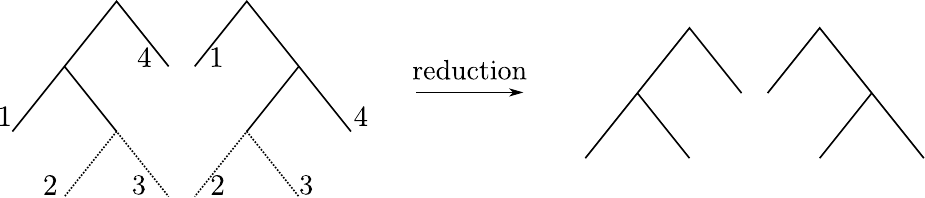}
	\end{center}
\caption{Example of reduction of carets}
\label{reduction_tree}
\end{figure}

Let $\mathcal{T}$ be the set of all pairs of binary trees whose number of leaves are the same. 
We define the equivalence relation on $\mathcal{T}$. 
Let $(T_+, T_-)$ be such two binary trees with $n$ leaves. 
We label the leaves of the trees with the numbers $1, 2, \dots, n$ from left to right, respectively. 
Assume that there exists $i \in \{1, 2, \dots, n-1\}$ such that $i$ and $i+1$ have a common parent in both $T_+$ and $T_-$. 
Then we can get the binary tree $T_+^\prime$ by removing leaves $i$, $i+1$ and corresponding two edges from $T_+$. 
Similarly, we get $T_-^\prime$ from $T_-$. 
See figure \ref{reduction_tree} for an example of this operation. 
We call this operation \textit{reduction of carets}. 
We say an element in $\mathcal{T}$ is \textit{reduced} if there exists no such $i$. 
Define the equivalence relation as the one generated by reductions and its inverses. 
It is known that there exists a unique reduced representative for each equivalence class \cite[\S 2]{cannon1996intro}. 

We call \textit{Thompson's group $F$} as the group of the quotient set of $\mathcal{T}$ equipped with the following product. 
Let $a=(A_+, A_-)$ and $b=(B_+, B_-)$ be in $\mathcal{T}$. 
By the previous operations, we get two element $(A_+^\prime, A_-^\prime)$ and $(B_+^\prime, B_-^\prime)$ which are equivalent to $a$ and $b$, respectively, and $A_-^\prime=B_+^\prime$ holds. 
Then we define the product of the equivalent class of $a$ and that of $b$ as that of $(A_+^\prime, B_-^\prime)$. 

The following fact is well known. 
\begin{proposition}[{\cite[\S 2]{cannon1996intro}}] 
Thompson's group $F$ is isomorphic to the group consisting of homeomorphisms on the closed interval $[0, 1]$ satisfying the following conditions: 
\begin{enumerate}
\item they are piecewise linear and preserve the orientation, 
\item in each linear part, its slope is a power of $2$, and
\item the breakpoints are in $\mathbb{Z}[\frac{1}{2}] \times \mathbb{Z}[\frac{1}{2}]$. 
\end{enumerate}
\begin{sproof}
Let $T$ be a binary tree. 
We decompose $[0, 1]$ by assigning a subinterval to each vertex of $T$. 
First, we consider the root to be $[0, 1]$. 
Next, if the parent is $[a, b]$, then we set the left child to $[a, (a+b)/2]$ and the right child to $[(a+b)/2, b]$, inductively. 
As a result, the leaves give the decomposition $[a_1, a_2], [a_2, a_3], \dots, [a_{n-1}, a_n]$ where $a_1=0$ and $a_n=1$. 
See also Figure \ref{tree_to_intervals}. 
\begin{figure}[tbp]
	\begin{center}
		\includegraphics[height=90pt]{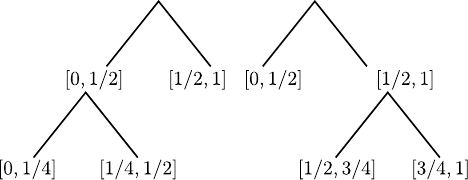}
	\end{center}
\caption{Vertices of a binary tree and subintervals}
\label{tree_to_intervals}
\end{figure}

Let $(T_+, T_-)$ be in $\mathcal{T}$. 
Since two trees have the same number of leaves, we get two decompositions $[a_1, a_2], [a_2, a_3], \dots, [a_{n-1}, a_n]$ and $[b_1, b_2], [b_2, b_3], \dots, [b_{n-1}, b_n]$. 
Therefore we get a piecewise linear map on $[0, 1]$ by mapping each $[a_i, a_{i+1}]$ linearly to $[b_i, b_{i+1}]$. 
This induces an isomorphism. 
\end{sproof}
\end{proposition}

We define virtual Thompson's group as a diagram group, but the following fact is helpful for understanding where its definition comes from. 
\begin{proposition}[{\cite[Example 6.4]{guba1997diagram}, \cite[Appendix]{farley2003finiteness}}]\label{F_iso}
Let $\cP_F=\langle x \mid x \to xx \rangle$ be a semigroup presentation. 
Then the diagram group $\cD(\cP_F, x)$ is isomorphic to the group $F$. 
\begin{sproof}
Let $\Delta$ be a reduced (geometric) diagram in $\cD(\cP_F, x)$. 
We construct a pair of binary trees from $\Delta$. 
This is achieved by associating each cell with a binary tree consisting of one parent and two children. 
Since each cell is of the form $x \to xx$ or $xx \to x$, we put a vertex on each edge and regard the vertex on the $x$-side (not $xx$-side) as the parent. 

By performing the same operation for all cells, we obtain a graph. 
Let $T_+$ and $T_-$ be the largest subgraphs whose roots are vertex on the top and bottom path, respectively. 
See Figure \ref{diagram_to_tree} and note that we omit the labels on all edges of diagrams since they are the same. 
\begin{figure}[tbp]
	\begin{center}
		\includegraphics[height=150pt]{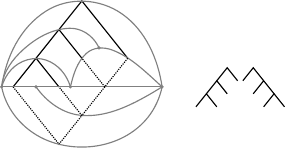}
	\end{center}
\caption{Example of the correspondence of a reduced diagram and a (reduced) pair of binary trees}
\label{diagram_to_tree}
\end{figure}
Since $\Delta$ is reduced, the union of $T_+$ and $T_-$ is the obtained graph, and their intersection is finitely many vertices. 
Moreover, $(T_+, T_-)$ is reduced. 
This induces an isomorphism. 
\end{sproof}
\end{proposition}
In the following, we define the virtual version of Thompson's group $F$. 
The name ``virtual'' comes from virtual knot theory described in Section \ref{subsection_virtualintro}. 
\begin{definition}
Let $\cP_{VF}$ be the semigroup defined by the following presentation: 
\begin{align*}
	\left \langle x, v \mathrel{}\middle | \mathrel{}
	\begin{array}{c}
		x \to xx, x \to vv, x \to vx, x \to xv \\ 
		v \to vv, v \to xx, v \to vx, v \to xv
	\end{array}
	\right \rangle. 
\end{align*}
Then we define \textit{virtual Thompson's group} $VF$ to be the diagram group $\cD(\cP_{VF}, x)$. 
\end{definition}

\begin{remark}
The group $VF$ is motivated by virtual knot theory as an analogy of the virtual braid group.
In this sense, this group is a  ``knot theoretic'' Thompson's group, and seems to be algebraically different from the so-called ``Thompson-like'' group.
\end{remark}

\subsection{Properties of $VF$}\label{subsection_properties_VF}
In this section, we list some properties of $VF$. 
The properties described below are already known to hold for diagram groups. 
See the respective references for details. 

The following statements $(1)$ and $(2)$ follow from \cite[Theorem 3.13, Theorem 4.1]{farley2003finiteness}, $(3)$ from \cite[Theorem 15.25]{guba1997diagram}, $(4)$ and $(5)$ from \cite[Theorem 6.1, Theorem 7.1]{guba2006diagram}, and (6) from \cite[Theorem 9.9]{guba2006diagram2}. 
\begin{theorem}\label{theorem_diagram_group}
Let $\cP$ be a semigroup presentation and $w$ be a given word. 
\begin{enumerate}[font=\normalfont]
\item If $\cP$ is a finite semigroup presentation, then $\cD(\cP, w)$ acts properly, cellularly, and freely by isometries on a proper $\mathrm{CAT(0)}$ cubical complex. 
\item If $\cP$ is a finite presentation and the semigroup is finite, then the diagram group $\cD(\cP, w)$ is of type $\mathscr{F}_{\infty}$. Especially, $\cD(\cP, w)$ is finitely presented. 
\item The group $\cD(\cP, w)$ has the unique extraction of root property. Especially, $\cD(\cP, w)$ is torsion-free. 
\item The group $\cD(\cP, w)$ is totally orderable. 
\item The group $\cD(\cP, w)$ is residually countable. 
\item All integer homology groups of $\cD(\cP, w)$ are free abelian. Especially, the abelianization of the group $\cD(\cP, w)$ is free abelian. 
\end{enumerate}
\end{theorem}
Note that statement $(1)$ has various corollaries such as satisfying Haagerup property and Baum--Connes conjecture. 
See \cite[Section 3.4]{farley2003finiteness} for details. 

Since $\cP_{VF}$ is a finite presentation of a finite semigroup, we have the following corollary: 
\begin{corollary}
The group $VF$ has all the properties in Theorem $\ref{theorem_diagram_group}$. 
\end{corollary}

In addition, by using only the relation $x \to xx$ in the rewriting, we have the following: 
\begin{proposition}
Thompson's group $F$ is a subgroup of $VF$. 
\end{proposition} 

In the rest of this section, we give an infinite presentation of $VF$ by using the Squier complex.
Let $\cP= \langle \Sigma \mid \cR \rangle$ be a semigroup presentation.
The \textit{Squier complex} $\cS(\cP)$ of $\cP$ is the 2-dimensional complex defined as follows:
\begin{itemize}
\item the 0-cells are the words on $\Sigma$;
\item the 1-cells $e_{p, u \to v, q}$ connect two 0-cells $puq$ and $pvq$ if $u \to v \in \cR$; and
\item the 2-cells $D_{p, u_1 \to v_1, q, u_2 \to v_2, r}$ bound the 4-cycles given by four 1-cells $e_{p, u_1 \to v_1, qu_2r}$, $e_{pv_1q, u_2 \to v_2, r}$, $e_{p, u_1 \to v_1, qv_2r}$, and $e_{pu_1q, u_2 \to v_2, r}$,
\end{itemize}
where $u_i \to v_i \in \cR\ (i=1,2)$ and $p, q, r$ are words on $\Sigma$.
For a given word $w$ on $\Sigma$, the diagram group $\cD(\cP, w)$ can be regarded as the fundamental group $\pi_1 (\cS(\cP), w)$.
See \cite[Section 6]{guba1997diagram} or \cite[Section 2]{genevois2022diagram} for details.

\begin{theorem} \label{presentation_VF}
The virtual Thompson's group $VF$ admits the following infinite presentation:
\begin{description}
\item[Generators:]\mbox{}
\begin{itemize}
\item $X_{v \to vv}, X_{v \to vx}, X_{v \to xv}$,
\item $X_{x,s \to t, u}, X_{v,s \to t, u}$\ \ $(s \in \{x, v\}, t \in \{xx, xv, vx, vv\}$ and $u$ is a word on $\Sigma)$,
\item $X_{x \to vv, u}, X_{x \to vx, u}, X_{v \to xx, u}, X_{v \to xv, u}$\ \ $(u$ is a non-empty word on $\Sigma)$.
\end{itemize}
\item[Relations:]\mbox{}
\begin{itemize}
\item $X_{x, s_1 \to t_1, ps_2q} X_{x, s_2 \to t_2, q} = X_{x, s_2 \to t_2,q} X_{x, s_1 \to t_1, pt_2q}$,
\item $X_{v, s_1 \to t_1, ps_2q} X_{v, s_2 \to t_2, q} = X_{v, s_2 \to t_2, q} X_{v, s_1 \to t_1, pt_2q}$,
\item $X_{x, s \to t, q} X_{x \to vv, ptq} = X_{x \to vv, psq} X_{v, s \to t, q}$,
\item $X_{x, s \to t, q} X_{x \to vx, ptq} = X_{x \to vx, psq} X_{v, s \to t, q}$,
\item $X_{v, s \to t, q} X_{v \to xx, ptq} = X_{v \to xx, psq} X_{x, s \to t, q}$,
\item $X_{v, s \to t, q} X_{v \to xv, ptq} = X_{v \to xv, psq} X_{x, s \to t, q}$,
\end{itemize}
where $s, s_1, s_2 \in\{x, v\}, t, t_1, t_2 \in \{xx, xv, vx, vv\}$, and $p, q$ are words in $\Sigma$.
\end{description}
\end{theorem}

For example, substituting $s_i = x, t_i = xx\ (i=1,2),\ p = x^j$, and $q = x^k\ (j, k \geq 0)$ in the first relation, we have
\begin{align*}
X_{x,x \to xx, x^j x x^k} X_{x, x \to xx, x^k} = X_{x, x \to xx, x^k} X_{x, x \to xx, x^j xx x^k}.
\end{align*}
Set $x_k \coloneqq X_{x, x \to xx, x^k}$, then we rewrite this relation as follows:
\begin{align*}
x_{j+k+1} x_k = x_k x_{j+k+2},
\end{align*}
that is,
\begin{align*}
x_n x_k = x_k x_{n+1}\ (0 \leq k < n),
\end{align*}
which is exactly the relation for Thompson's group $F$.

From the presentation in Theorem \ref{presentation_VF}, three generators $X_{v \to vv}, X_{v \to vx}$ and $X_{v \to xv}$ have no relations.
Therefore, the virtual Thompson's group $VF$ is the free product of the free group of rank 3 generated by these generators and the remaining part of $VF$.

\begin{figure}[tbp]
\begin{center}
\includegraphics[width=150mm]{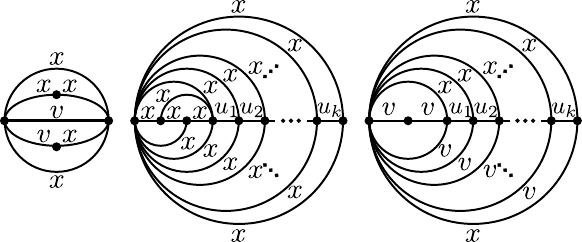}
\end{center}
\caption{The generators $X_{x \to vx}, X_{x, x\to xx, u}$ and $X_{z \to vv, u}$ where $u$ is a word $u_1 u_2 \cdots u_k$. }
\label{generators_VF}
\end{figure}

Finally, the generators of $VF$ can be described as the geometric diagrams shown in Figure \ref{generators_VF}.

\section{The Construction of virtual links from virtual Thompson's group}\label{section_construction}
\subsection{The construction} \label{construction}
In this section, we explain the construction of a virtual link from an element of virtual Thompson's group $VF$ with an example.
This construction is based on \cite{jones2017thompson} and \cite{golan2017jones}.

\noindent
\textbf{Step 1: Construct the Thompson graph $T(\Delta)$.}

Let $\Delta$ be a reduced diagram in $VF = \cD(\cP_{VF},x)$.
We define the \textit{Thompson graph} $T(\Delta)$ as a ``subgraph'' of the diagram $\Delta$ as follows (cf. \cite[Definition 3.2]{golan2017jones}):
the vertices of $T(\Delta)$ are all vertices of $\Delta$ except the terminal vertex.
In order to define the edges of $T(\Delta)$, we use the following lemma.

\begin{figure}[tbp]
\begin{center}
\includegraphics[height=60pt]{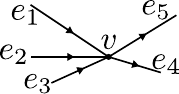}
\end{center}
\caption{The edge $e_1$ is the first incoming edge and $e_3$ is the last incoming edge of $v$.}
\label{incoming_outcoming}
\end{figure}

\begin{lemma}[{\cite[Lemma 3.7]{guba1997diagram}}] \label{order}
For any inner vertex $v$ of $\Delta$, that is, the vertex which does not coincide with the initial vertex nor the terminal vertex, there uniquely exists a sequence $e_1, \ldots, e_n$ of edges with endpoint $v$ in the counterclockwise order such that for some $k\ (1 \leq k < n)$, edges $e_1, \ldots, e_k$ are incoming and edges $e_{k+1}, \ldots, e_n$ are outgoing $($see Figure $\ref{incoming_outcoming})$.
\end{lemma}

\begin{figure}[tbp]
	\begin{center}
		\includegraphics[height=150pt]{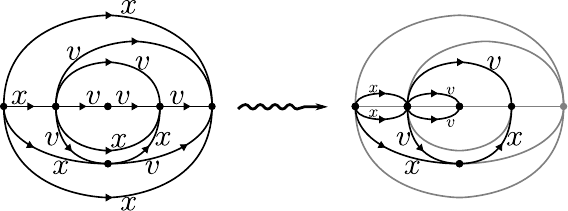}
	\end{center}
\caption{An example of the Thompson graph $T(\Delta)$}
\label{example_step1}
\end{figure}

For any inner vertex $v$ of $\Delta$, we assign numbers to edges with endpoint $v$ as in Lemma \ref{order}.
The edges of $T(\Delta)$ are the first and the last incoming edges with respect to the order for each inner vertex of $\Delta$.
If the first and the last edges of $v$ coincide, that is, $v$ has exactly one incoming edge, then we make a copy of the incoming edge labeled by the same letter.
Therefore, any inner vertex of $T(\Delta)$ has two incoming edges.
Figure \ref{example_step1} is an example of this step.

\noindent
\textbf{Step 2: Construct the medial graph $M(T(\Delta))$.}

\begin{figure}[tbp]
\begin{center}
\includegraphics[height=110pt]{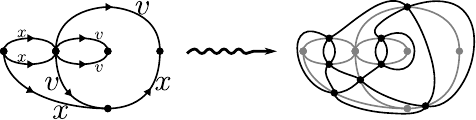}
\end{center}
\caption{An example of the medial graph $M(T(\Delta))$}
\label{example_step2}
\end{figure}

The medial graph is defined for any connected plane graph.
Let $G$ be a connected plane graph, and then its medial graph $M(G)$ is obtained as follows:
We put a vertex of $M(G)$ on every edge of $G$, and join two new vertices by an edge if the corresponding edges of $G$ are adjacent on a face of $G$.
Figure \ref{example_step2} is an example of this step.

\noindent
\textbf{Step 3: Construct the virtual link diagram $L(\Delta)$.}

\begin{figure}[tbp]
	\begin{center}
		\includegraphics[height=110pt]{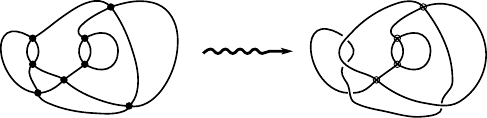}
	\end{center}
\caption{An example of the virtual link diagram $L(\Delta)$}
\label{example_step3}
\end{figure}

In general, because the medial graph is 4-valent, we are able to obtain a virtual link diagram $L(\Delta)$ by turning all vertices of $M(T(\Delta))$ into classical or virtual crossings:
if the corresponding edge in $T(\Delta)$ is
\begin{align} \label{crossing_rule}
\begin{cases}
\textrm{the first and labeled by $x$, then \cinput{3}{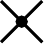} $\to$ \cinput{3}{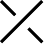},} \\
\textrm{the last and labeled by $x$, then \cinput{3}{vertex.pdf_tex} $\to$ \cinput{3}{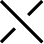}, and} \\
\textrm{labeled by $v$, then \cinput{3}{vertex.pdf_tex} $\to$ \cinput{3}{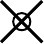}.}
\end{cases}
\end{align}
Figure \ref{example_step3} is an example of this step. 

\subsection{Labeled binary trees} \label{labeled}
\begin{figure}[tbp]
	\begin{center}
		\includegraphics[height=140pt]{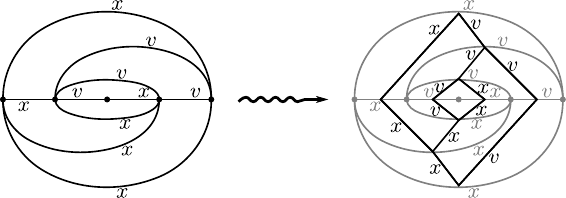}
	\end{center}
\caption{An example of the correspondence of a diagram and a pair of labeled binary trees. }
\label{VFdiagram_to_VFtree}
\end{figure}

In this section, we discuss the relationship between elements of $VF$ and labeled binary trees.
Suppose that the diagram $\Delta: x = w_1 \to w_2 \to \cdots \to w_{n-1} \to w_n = x$ satisfies the following condition:
\begin{equation} \label{condition}
\begin{split}
&\text{There uniquely exists $i \in \mathbb{Z} \cap [1, n]$ such that}
\begin{cases}
|w_j| < |w_{j+1}| & (1 \leq j <i) \\
|w_j| > |w_{j+1}| & (i \leq j < n)
\end{cases} \text{hold},
\end{split}
\end{equation}
where $| \cdot |$ denotes the length of a word.
Geometrically, this condition implies that all vertices of $\Delta$ can be placed on a straight line, and every vertex except the initial vertex has an incoming edge connected to its immediate left one.
The path of $\Delta$ on the straight line is exactly the trivial geometric diagram of $w_i$.
Then the cell $w_j \to w_{j+1}$ for $1 \leq j < i$ is the $(s,t)$-cell and the cell $w_j \to w_{j+1}$ for $i \leq j < n$ is the $(t,s)$-cell, where $s \in \{x, v\}$ and $t \in \{xx, xv, vx, vv\}$.
In this case, the first and last edges coincide with the top-most and bottom-most incoming edges of \cite[Definition 3.2]{golan2017jones}, respectively.
Similarly to the proof of Proposition \ref{F_iso}, the diagram $\Delta$ can be described as a pair $(T_+, T_-)$ of labeled binary trees with the same number of leaves.
The label of each edge is determined by the label of the corresponding ``child'' edge of the cell in $\Delta$.
We give an example in Figure \ref{VFdiagram_to_VFtree}.

\begin{figure}[tbp]
	\begin{center}
		\includegraphics[height=140pt]{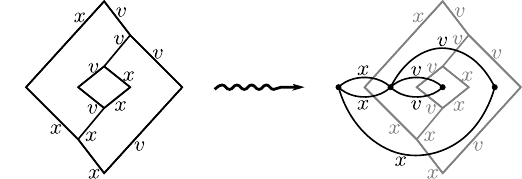}
	\end{center}
\caption{An example of the correspondence of a pair of labeled binary trees and a $\Gamma$-graph. }
\label{VFtree_to_Jones}
\end{figure}

On the other hand, Jones \cite{jones2017thompson} introduced a method of constructing a link diagram from an element of $F$ by using a pair of binary trees.
In the case above, this construction can be extended naturally.
Let $(T_+, T_-)$ be a pair of reduced labeled binary trees with $n+1$ leaves obtained from an element of $VF$, and place its leaves at $\left(\frac{1}{2}, 0 \right), \left(\frac{3}{2}, 0 \right), \ldots, \left(\frac{2n+1}{2}, 0 \right)$.
Note that the tree $T_+$ is in the upper half-plane, and $T_-$ is in the lower half-plane.
The plane graph $\Gamma(T_+, T_-)$, which is called the $\Gamma$-graph of $(T_+, T_-)$, is defined as follows: the vertices of $\Gamma(T_+, T_-)$ are put at $(0, 0), (1, 0), \ldots, (n, 0)$.
An edge of $\Gamma(T_+, T_-)$ passes transversely just once an edge \cinput{2}{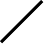} of $T_+$ (i.e. an edge from top right to bottom left) or an edge \cinput{3}{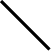} of $T_-$ (i.e. an edge from top left to bottom right) and does not do the other edges of  $(T_+, T_-)$.
Every edge is labeled by $x$ or $v$ corresponding to the label of an edge of $(T_+, T_-)$.
We illustrate an example in Figure \ref{VFtree_to_Jones}. 

For the two constructions, the following holds: 
\begin{proposition}[cf. {\cite[Proposition 3.5]{golan2017jones}}] \label{lemma_Thompsongraph_Jones}
Let $\Delta$ be a diagram of $VF$ satisfying condition $(\ref{condition})$ and $(T_+, T_-)$ the associated graph.
Then the Thompson graph $T(\Delta)$ is isomorphic to the $\Gamma$-graph $\Gamma(T_+, T_-)$.
\end{proposition}
\begin{sproof}
Golan and Sapir \cite{golan2017jones} proved this lemma in the case of $F$ by stretching the edges of $T(\Delta)$.
From their argument, we see that if the diagram $\Delta$ of $F$ is not reduced but satisfies condition (\ref{condition}), then $T(\Delta)$ is isomorphic to $\Gamma(T_+, T_-)$.
In the case of $VF$, the correspondence of the labels of the edges of $T(\Delta)$ and $\Gamma(T_+, T_-)$ is clear from the construction. 
We illustrate an example in Figure \ref{example_lemma3.2}. 
\begin{figure}[tbp]
	\begin{center}
		\includegraphics[height=140pt]{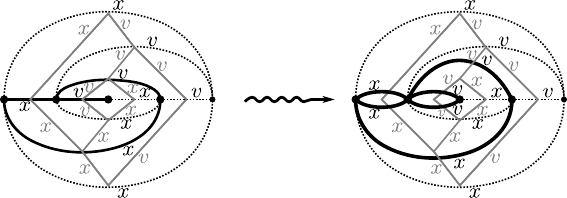}
	\end{center}
\caption{An example of the correspondence in Lemma \ref{lemma_Thompsongraph_Jones}. By stretching the edges of $T(\Delta)$, it is isomorphic to the graph $\Gamma(T_+, T_-)$. The gray letters are labels of binary trees. }
\label{example_lemma3.2}
\end{figure}
\end{sproof}

\section{Proof of Theorem \ref{main_theorem}} \label{section_proof_main_theorem}
Theorem \ref{main_theorem} states that every virtual link can be described as the virtual link diagram $L(\Delta)$ for an element $\Delta$ of $VF$.
In this section, we prove this theorem.
The procedure of the proof is based on \cite{jones2017thompson}.
In fact, we are always able to choose an element of $VF$ representing a given virtual link which satisfies condition $(\ref{condition})$.

Let $L$ be a virtual link diagram.

\noindent
\textbf{Step 1: Construct the Tait graph $T(L)$.}

\begin{figure}[tbp]
\begin{center}
\includegraphics[height=70pt]{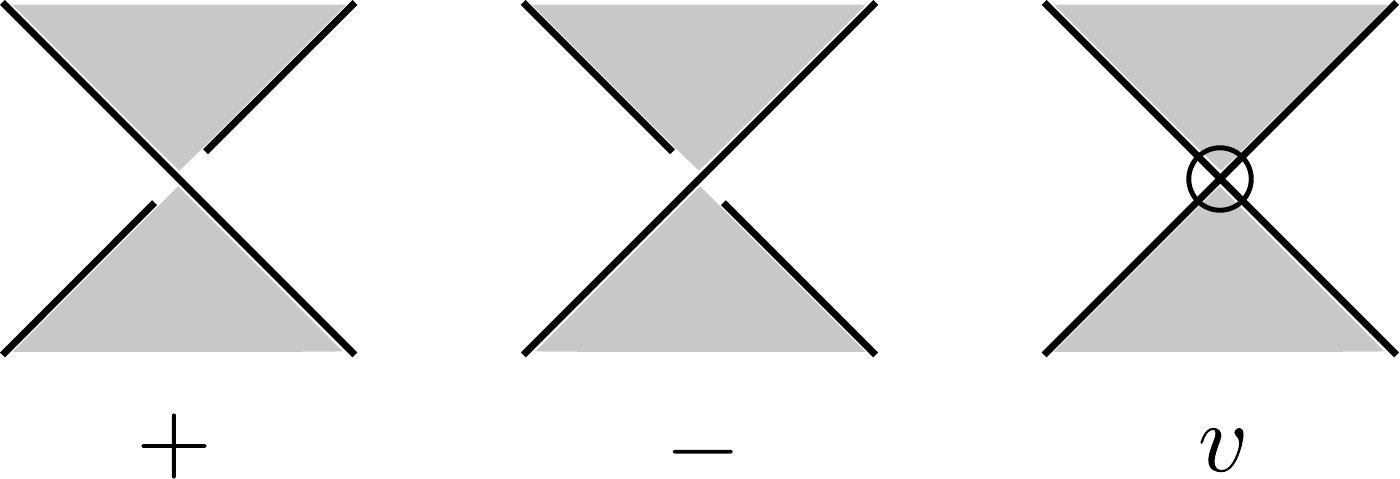}
\end{center}
\caption{The labels of crossings}
\label{checkerboard_crossing}
\end{figure}

We apply the checkerboard coloring to the diagram $L$, that is, we paint regions of $L$ with black or white so that adjacent regions are different colors.
By convention, the color of the unbounded region is white.
The vertices of the Tait graph $T(L)$ correspond to the black regions of $L$, and the edges correspond to the crossings and are labeled by $+$, $-$, or $v$ according to the rule in Figure \ref{checkerboard_crossing}.

\begin{figure}[tbp]
\begin{center}
\includegraphics[height=110pt]{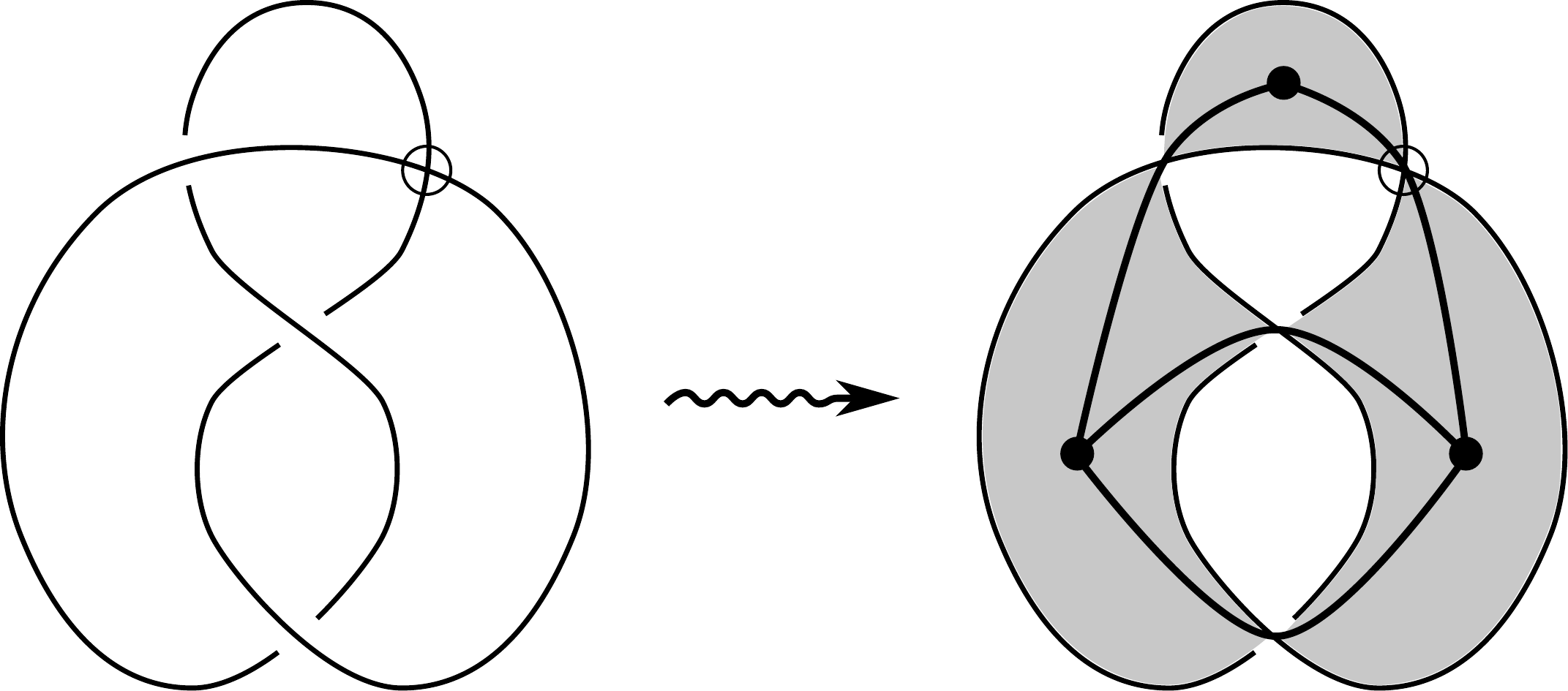}
\end{center}
\caption{An example of the checkerboard coloring and the Tait graph}
\label{checkerboard_example}
\end{figure}

\noindent
\textbf{Step 2: Deform the graph $T(L)$.}

Jones \cite{jones2017thompson} showed the condition that a connected plane graph can be obtained from an element $(T_+, T_-)$ of $F$:

\begin{lemma}[{\cite[Lemma 4.1.4]{jones2017thompson}}] \label{property}
Let $\Gamma$ be a connected plane graph.
Suppose that $\Gamma$ consists of two trees, $\Gamma_+$ in the upper half-plane and $\Gamma_-$ in the lower half-plane, and these two trees satisfy the following properties:
\begin{enumerate}
\item the vertices are $(0,0), (1,0), \ldots, (n,0)$, \label{property1}
\item each vertex other than $(0,0)$ is connected to exactly one vertex to its left one, and \label{property2}
\item each edge can be parametrized by a smooth curve $(x(t), y(t))$ for $t \in [0,1]$ with $x'(t)>0$ and either $y(t)>0$ or $y(t)<0$ for $t \in (0,1)$. \label{property3}
\end{enumerate}
Then there exists an element $(T_+, T_-)$ of $F$ such that $\Gamma(T_+, T_-)$ is isomorphic to $\Gamma$.
Equivalently, there exists a diagram $\Delta$ of $F$ such that $T(\Delta)$ is isomorphic to $\Gamma$.
\end{lemma}

For the virtual case, the Thompson graph has labels $x$ or $v$.
In particular, if there exists a vertex of a diagram $\Delta$ with exactly one incoming edge, then it has two incoming edges with the same labels in the Thompson graph $T(\Delta)$.
Hence, we obtain the condition for the virtual version:

\begin{lemma} \label{virtual_property}
Let $\Gamma$ be a connected plane graph with each edge labeled by $x$ or $v$.
Suppose that $\Gamma$ consists of two trees, $\Gamma_+$ in the upper half-plane and $\Gamma_-$ in the lower half-plane, and these two trees satisfy the properties $(\ref{property1})$, $(\ref{property2})$ and $(\ref{property3})$ in Lemma $\ref{property}$.
Moreover, assume that $\Gamma$ satisfies the following condition:
\begin{enumerate}
\setcounter{enumi}{3}
\item Two edges connecting adjacent two vertices have the same labels. \label{property4}
\end{enumerate}
Then there exists a pair $(T_+, T_-)$ of labeled binary trees in $VF$ such that $\Gamma(T_+, T_-)$ is isomorphic to $\Gamma$.
Equivalently, there exists an element $\Delta$ of $VF$ satisfying condition $(\ref{condition})$ such that $T(\Delta)$ is isomorphic to $\Gamma$.
\end{lemma}

For a diagram $\Delta$ satisfying condition (\ref{condition}), the Tait graph of $L(\Delta)$ (i.e. the Thompson graph $T(\Delta)$) satisfies the conditions of Lemma \ref{virtual_property}, with edges in the upper half-plane labeled by $+$ or $v$ and edges in the lower half-plane labeled by $-$ or $v$.
Therefore, in order to prove the main theorem, we apply the Reidemeister moves on the given Tait graph so that the deformed graph is of the form we have just described.

\begin{figure}[tbp]
\begin{center}
\includegraphics[height=180pt]{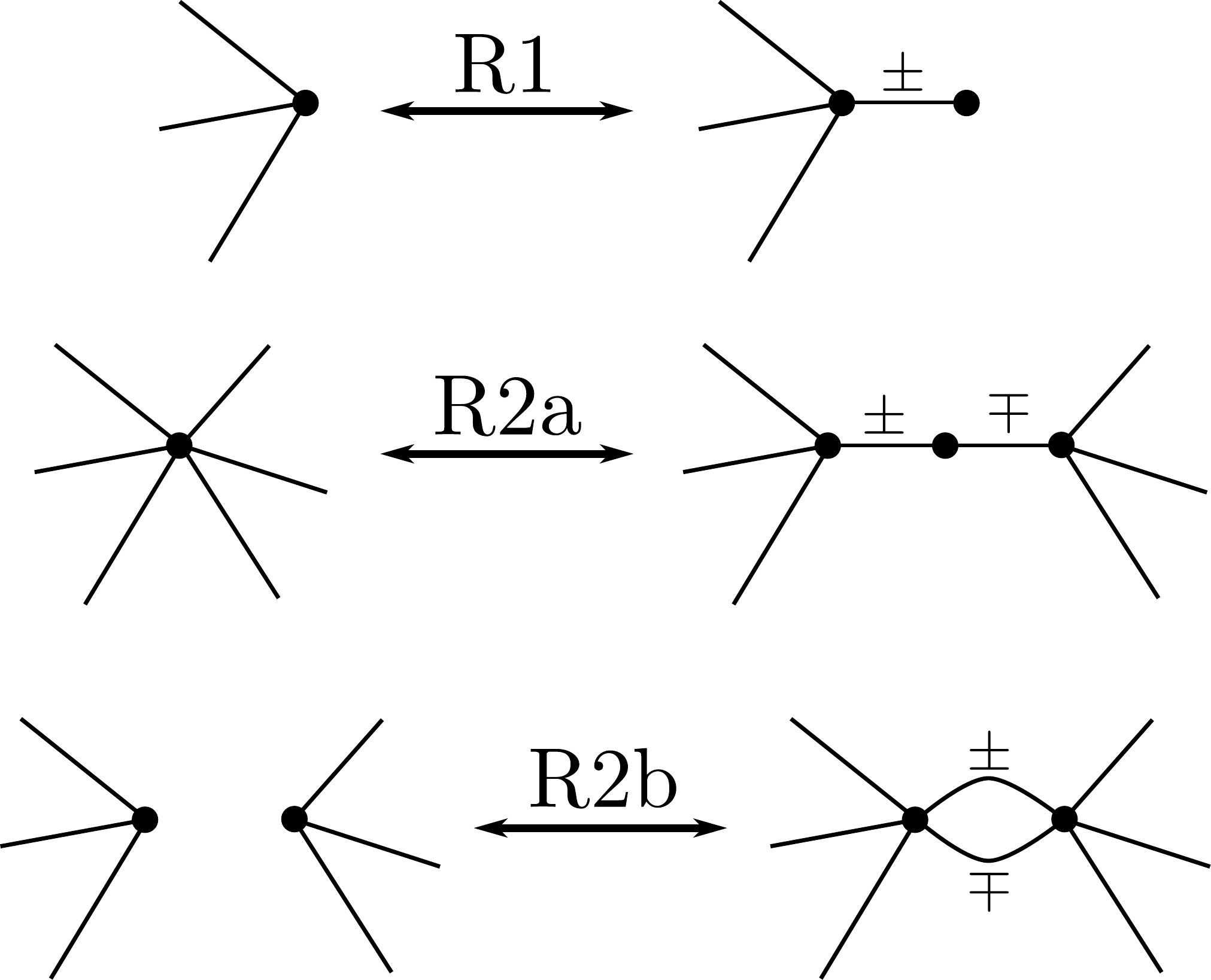}
\end{center}
\caption{The Reidemeister moves R1 and R2 on the plane graph}
\label{graph_Reidemeister}
\end{figure}

We recall some local moves on the labeled plane graph corresponding to the Reidemeister moves R1 and R2 (see Figure \ref{graph_Reidemeister}).

\begin{definition}[{\cite[Definition 5.3.4]{jones2017thompson}}]
Two labeled plane graphs are \textit{2-equivalent} if they differ by planar isotopies and any of the moves R1, R2a, and R2b.
\end{definition}

\begin{lemma}[{\cite[Lemma 5.3.6]{jones2017thompson}}] \label{standard}
Any Tait graph is $2$-equivalent to a plane graph satisfying conditions $(\ref{property1})$ and $(\ref{property3})$ in Lemma $\ref{property}$.
\end{lemma}

Therefore, we may assume that the Tait graph $T(L)$ satisfies the above conditions.
Suppose that the edges of a standard plane graph are oriented from left to right.
We recall some notations in \cite{jones2017thompson}.

\begin{definition}[{\cite[Definition 5.3.7 and 5.3.8]{jones2017thompson}}]
For a vertex $u$ of $T(L)$, we set
\begin{align*}
e^{\textrm{up}} &:= \left\{ e \in E(T(L)) \mathrel{}\middle | \mathrel{} e\ \textrm{lies in the upper half-plane} \right\},\\
e^{\textrm{down}} &:= \left\{ e \in E(T(L)) \mathrel{}\middle | \mathrel{} e\ \textrm{lies in the lower half-plane} \right\},\\
e^{\textrm{in}}_u &:= \left\{ e \in E(T(L)) \mathrel{}\middle | \mathrel{} \tau(e)=u \right\},\\
e^{\textrm{out}}_u &:= \left\{ e \in E(T(L)) \mathrel{}\middle | \mathrel{} \iota(e)=u \right\},
\end{align*}
where $E(T(L))$ is the set of edges of $T(L)$, and $\tau(e)$ and $\iota(e)$ are the terminal and initial vertices of $e$, respectively.
\end{definition}

\noindent
\textbf{Case 1.}
There exists a vertex $u$ different from $(0,0)$ with $e^{\textrm{in}}_u = \emptyset$.
Let $w$ be the vertex immediately to the left of $u$ as below:
\begin{align*}
\cinput{0}{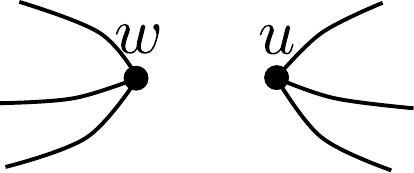}
\end{align*}
We add two edges connecting $w$ and $u$ so that the deformed graph is 2-equivalent to the original graph:
\begin{align*}
\cinput{0}{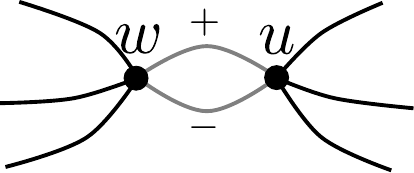}
\end{align*}

\noindent
\textbf{Case 2.}
There exists a vertex $u$ with $|e^{\textrm{in}}_u| = 1$.
We may assume that the incoming edge of $u$ is in the upper half-plane.
This is labeled by $+$, $-$, or $v$.
The situation near $u$ is as below:
\begin{align*}
\cinput{0}{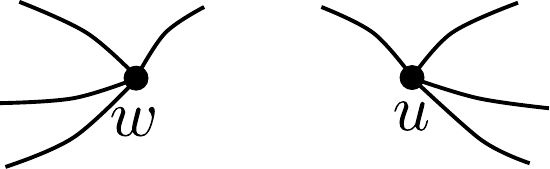}
\end{align*}
Then we add one vertex and three edges as below:
\begin{align*}
\cinput{0}{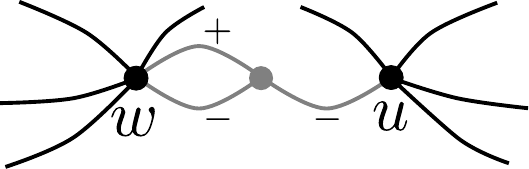}
\end{align*}

\noindent
\textbf{Case 3.}
There exists a vertex $u$ with $|e^{\textrm{in}}_u \cap e^{\textrm{up}}| > 1$ or $|e^{\textrm{in}}_u \cap e^{\textrm{down}}| > 1$.
We may show only the first case, and the other case is similar.
The situation near $u$ is below:
\begin{align*}
\cinput{0}{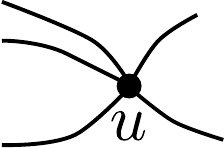}
\end{align*}
Then we add three vertices and five edges as below:
\begin{align*}
\cinput{0}{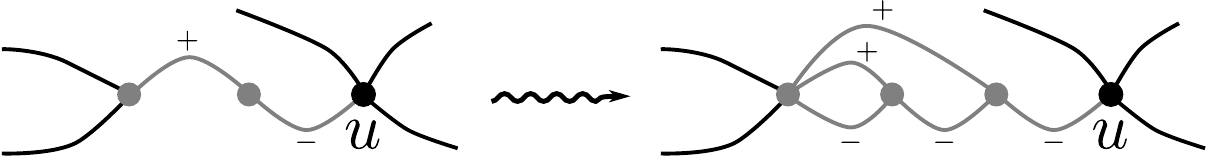}
\end{align*}

\noindent
\textbf{Case 4.}
After applying the previous deformations, all vertices, except the vertex $(0,0)$, have two incoming edges, one in the upper half-plane and the other in the lower half-plane.
Hence, this graph satisfies condition (\ref{property2}) in Lemma \ref{property}.
Then we may have two problems that
\begin{enumerate}
\item[(i)] a $-$-labeled edge is in the upper half-plane or a $+$-labeled edge is in the lower half-plane, and
\item[(ii)] two edges connecting the adjacent two vertices have labels $+$ and $v$, or $v$ and $-$, respectively.
\end{enumerate}
We consider the first problem, and we may show only the case of $-$-labeled edge in the upper half-plane.
This situation looks like:
\begin{align*}
\cinput{0}{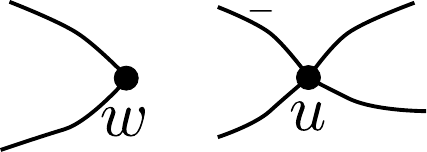}
\end{align*}
Then we apply the following deformation
\begin{align*}
\cinput{0}{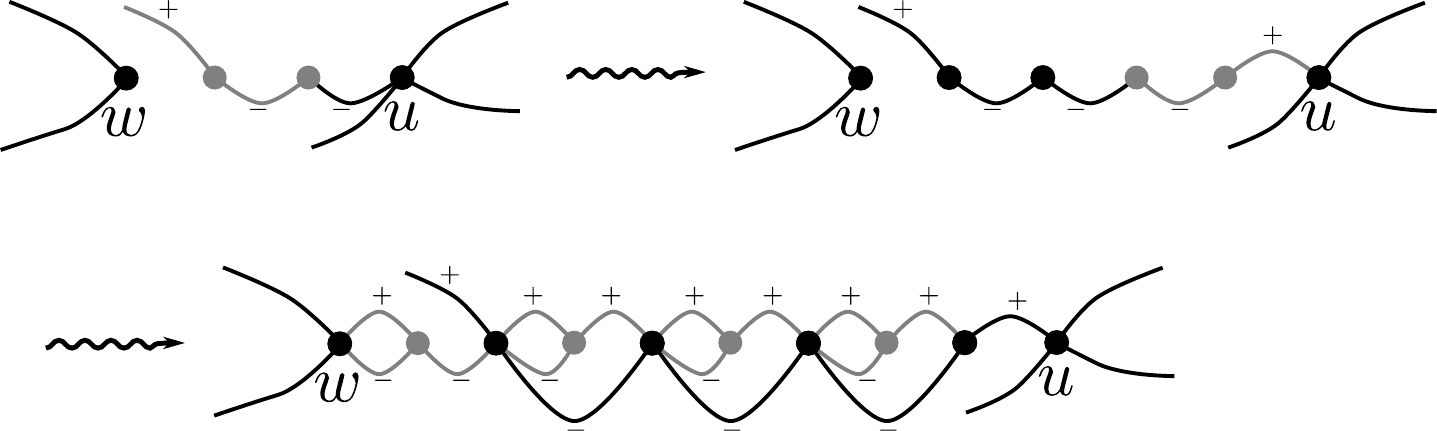}
\end{align*}

Next, we consider the second problem.
We may show only the case that two labels are $+$ and $v$, respectively.
Suppose that such two edges connect the adjacent vertices $w$ and $u$, then the situation looks like:
\begin{align*}
\cinput{0}{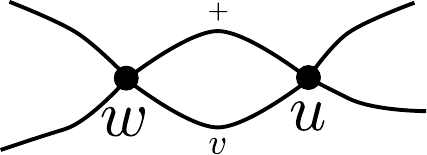}
\end{align*}
Then we apply the following deformations:
\begin{align*}
\cinput{0}{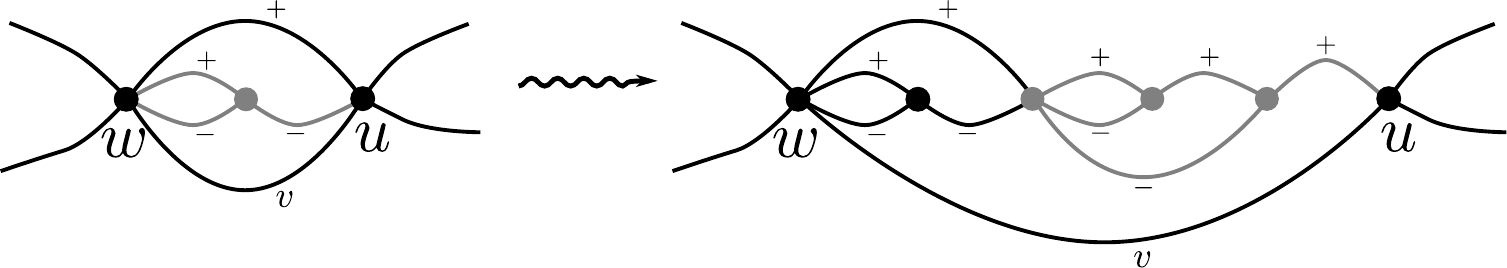}
\end{align*}

From the above, the proof of Theorem \ref{main_theorem} is complete.

\begin{example}
Figure \ref{algorithm_example} shows the application of the algorithm to the virtual knot $3.1$ in the list\footnote{\url{http://www.math.toronto.edu/~drorbn/Students/GreenJ/}} by Jeremy Green.
Its last figure is a diagram of $VF$ representing $3.1$.
The top and bottom paths must be labeled by $x$ from the definition.
However, other than those edges, gray edges can be labeled by either $x$ or $v$.
\begin{figure}[tbp]
\begin{center}
\includegraphics{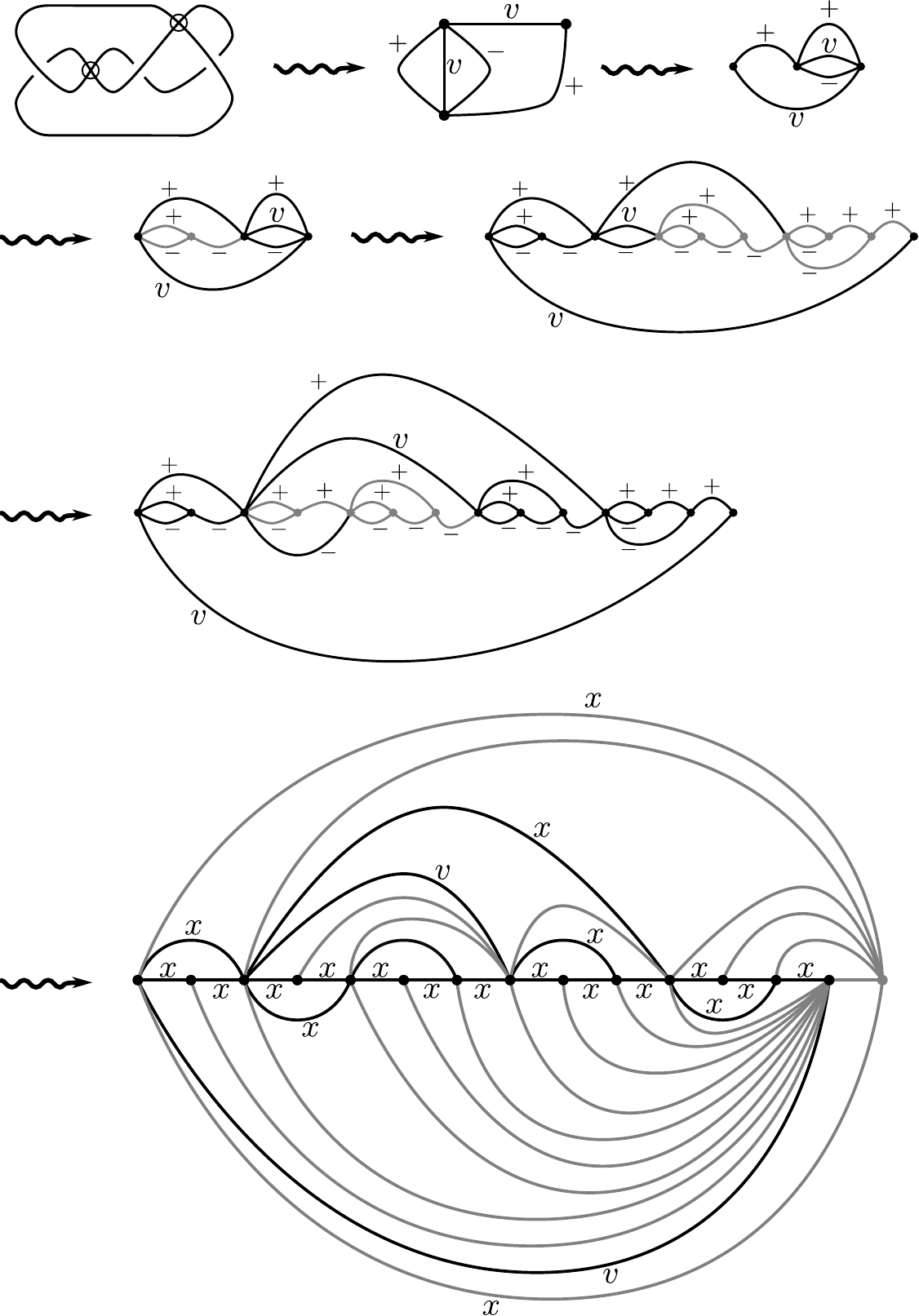}
\end{center}
\caption{}
\label{algorithm_example}
\end{figure}
\end{example}

Since a virtual link is an immersion of circles, its orientation is induced from the one of each circle.
Jones defined a subgroup $\overrightarrow{F}$ of $F$ which is called oriented Thompson's group.
This group consists of all pairs of binary trees whose $\Gamma$-graphs are 2-colorable, and its element yields an oriented link.
Aiello \cite{aiello2020alexander} proved Alexander's theorem for the oriented case by using another local move.
From \cite[Lemma 4.1]{golan2017jones}, we are able to define a subgroup $\overrightarrow{VF}$ of $VF$ consisting of all diagrams whose Thompson graphs are 2-colorable.
Moreover, by using Aiello's move, the oriented version of Theorem \ref{main_theorem} can be proved similarly.

\begin{theorem}\label{main_theorem_oriented}
Any oriented virtual link can be obtained from an element in $\overrightarrow{VF}$.
\end{theorem}

Finally, Golan and Sapir \cite{golan2017jones} showed that oriented Thompson's group $\overrightarrow{F}$ is isomorphic to 3-adic Thompson's group $F(3)$, which is a diagram group especially.
In general, a subgroup of the diagram group is not always a diagram group, and thus there is a natural problem whether $\overrightarrow{VF}$ is a diagram group or not.

\section*{Acknowledgements}
We would like to thank Professor Tomohiro Fukaya who is the first author's supervisor for his several comments. 
We also wish to thank Professor Takuya Sakasai who is the second author's supervisor and Xiaobing Sheng Ph.\ D.\ for their helpful comments.
We are also very grateful to Anthony Genevois for the valuable comments.

\bibliographystyle{plain}
\bibliography{bib1} 

\bigskip
\address{
DEPARTMENT OF MATHEMATICAL SCIENCES,
TOKYO METROPOLITAN UNIVERSITY,
MINAMI-OSAWA HACHIOJI, TOKYO, 192-0397, JAPAN
}

\textit{E-mail address}: \href{mailto:kodama-yuya@ed.tmu.ac.jp}{\texttt{kodama-yuya@ed.tmu.ac.jp}}

\address{GRADUATE SCHOOL OF MATHEMATICAL SCIENCES, THE
UNIVERSITY OF TOKYO, 3-8-1 KOMABA, MEGURO-KU, TOKYO, 153-8914,
JAPAN}

\textit{E-mail address}: \href{takano@ms.u-tokyo.ac.jp}{\texttt{takano@ms.u-tokyo.ac.jp}}
\end{document}